# A BI-OBJECTIVE STOCHASTIC APPROACH FOR THE STOCHASTIC CARP

**Gérard Fleury, Philippe Lacomme, Christian Prins, Marc Sevaux**

# A Bi-Objective Stochastic Approach for the Stochastic CARP

GERARD FLEURY
*Laboratoire de Mathématiques, UMR 6620*
*Université Blaise Pascal, Campus des Cézeaux*
*63177 Aubière Cedex*

PHILIPPE LACOMME
*Laboratoire d'Informatique (LIMOS), UMR 6158*
*Université Blaise Pascal, Campus des Cézeaux*
*63177 Aubière Cedex*

CHRISTIAN PRINS
*Institut Charles Delaunay*
*Université de Technologie de Troyes (UTT), BP 2060*
*10010 Troyes Cedex, France*

MARC SEVAUX
*Université de Bretagne Sud*
*CNRS, FRE 2734 – LESTER*
*Centre de recherche - BP 92116*
*56321 LORIENT Cedex France*

Abstract

---

The Capacitated Arc Routing Problem (CARP) occurs in applications like urban waste collection or winter gritting. It is usually defined in literature on an undirected graph $G=(V,E)$, with a set $V$ of $n$ nodes and a set $E$ of $m$ edges. A fleet of identical vehicles of capacity $Q$ is based at a depot node. Each edge $i$ has a cost (length) $c_i$ and a demand $q_i$ (*e.g.* an amount of waste), and it may be traversed any number of times. The edges with non-zero demands or tasks require service by a vehicle. The goal is to determine a set of vehicle trips (routes) of minimum total cost, such that each trip starts and ends at the depot, each task is serviced by one single trip, and the total demand handled by any vehicle does not exceed $Q$. To the best of our knowledge the best published method is a memetic algorithm first introduced in 2001.
This article provides a new extension of the NSGA II (Non-dominated Sorting Genetic Algorithm) template to comply with the stochastic sight of the CARP. The main contribution is:

- to introduce mathematical expression to evaluate both cost and duration of the longest trip and also standard deviation of these two criteria.
- to use a NGA-II template to optimize simultaneously the cost and the duration of the longest trip including standard deviation.

The numerical experiments managed on the thee well-known benchmark sets of DeArmon, Belenguer and Benavent and Eglese, prove it is possible to obtain robust solutions in four simultaneous criteria in rather short computation times.

Key Words: CARP, bi-objective, stochastic

---

## 1. INTRODUCTION

### 1.1. The basic CARP

The *Capacitated Arc Routing Problem* (CARP) consists of visiting a subset of edges instead of the nodes as the well-known VRP. CARP applications include for instance urban waste collection, winter gritting and inspection of power lines. To make the paper more concrete, and without loss of generality, examples are inspired by urban waste collection.

The CARP of literature tackles undirected networks. Each edge models a two-way street which both sides are treated in parallel and in any direction (*bilateral collection*), a common practice in residential areas with narrow streets. A fleet of identical vehicles of limited capacity is based at a depot node. Each edge can be traversed any times, with a known traversal cost. Some edges are required, *i.e.*, they have a non-zero demand (amount of waste) to be collected by a vehicle. The CARP consists in determining a set of vehicle trips minimizing the total cost, such that each trip starts and ends at the depot node. Each required edge is serviced by one single trip, and therefore must be visited at most once. The total demand processed by a trip must not exceed the vehicle capacity. Each vehicle has the same capacity.

The CARP is NP-hard, even in the single-vehicle case called Rural Postman Problem (RPP). Since exact methods are still limited to 20-30 edges (Hirabayashi *et al.*, 1992), heuristics are required for large instances, *e.g.* Augment-Merge (Golden and Wong, 1981), Path-Scanning (Golden *et al.*, 1983), Construct-and-strike (Pearn's improved version, 1989), Augment-Insert (Pearn, 1991) and Ulusoy's tour splitting algorithm (Ulusoy, 1985). The first metaheuristic for the CARP, a simulated annealing procedure, was designed by Eglese in 1994 for winter gritting problems. Several tabu search (TS) algorithms are also available, both for special cases like the undirected RPP (Hertz *et al.*, 1999) or the mixed RPP (Corberan *et al.*, 2000) and for the CARP itself (Eglese and Li, 1996) (Hertz *et al.*, 2000). All these metaheuristics and classical heuristics can be evaluated thanks to lower bounds, generally based on linear programming formulations (Benavent *et al.*, 1992), (Belenguer and Benavent, 2003), (Amberg and Voβ, 2002). On most instances, the best-known lower bound is obtained by a cutting-plane algorithm (Belenguer and Benavent, 2003). To the best of our knowledge the best previous published method seems to be the memetic algorithm first introduced in 2001 by (Lacomme *et al.*, 2001) which outperforms the well-known CARPET method of Hertz.

### 1.2. The Stochastic CARP

The SCARP problem is similar to a CARP problem, except that positive demands $q_{ij}$ then become positive random variables $Q_{ij}$. To any SCARP problem, can be canonically associated a CARP problem, where the stochastic demands $Q_{ij}$ are replaced by their expectation $q_{ij} = \overline{Q}_{ij}$. To avoid any ambiguity we call "stochastic" an element applied to the SCARP and "deterministic" an element applied to its associated CARP. The objective in the CARP consists in determining a set of trips of minimal cost, and the objective solving the SCARP consists in determining a robust set of trips (a robust solution). Based on the definition of Jensen (Jensen, 2001), robust solutions are solutions well performing in front of variations in quantities to collect. Let us note:

$Q$ the common capacity of the vehicles used.

$c_i$ the (deterministic) cost of the edge $i$ ($c_i > 0$).

$q_i$ the deterministic demand on edge $i$ ($0 < q_i \le Q$).

$Q_i$ the random demand on edge $i$ ($0 < Q_i \le Q$).

$x$ a (deterministic) solution of the CARP (finite set of trips).

| | |
|---|---|
| $X$ | the finite set of all solutions. |
| $t$(x) | the deterministic number of trips of $x$ . |
| $T(x,\omega)$ | the random number of trips of $x$ (depending on random variation $\omega$ of demands). |
| $h(x)$ | the deterministic cost of a solution $x$ . |
| $H(x,\omega)$ | the random variable being the cost of $x$ (depending on random variation $\omega$ of demands). The cost is the sum of the costs of every edge used by each vehicle, with, sometimes, extra trips to the depot node... |

Let us consider a (deterministic) solution. If the random demands serviced by a trip become less important than expected, the cost does not vary, but, if the vehicle exhausts its capacity before the end of the trip, it must move from its current position in the network to the depot node and turn back to complete the trip initially planned. Such operations create an extra trip and imply a (possibly huge) increasing of the total cost. Such events occur in many applications like waste collection. Calling another vehicle may be impossible for several reasons: the driver can not inform his colleagues due to the lack of communication systems, the other drivers can not come because all trips are performed in parallel, the driver is the only one knowing this sector, and so on. Therefore, for any solution $x$ , when $q_i = \overline{Q_i}$ , one has $H(x,\omega) \geq h(x)$ .

Let us consider, for example, a CARP instance with only 9 tasks, 3 vehicles with $Q = 4$ and demands equal to 1 for each task. Assume a solution with a deterministic cost $h(x) = 100$ and 3 trips (figure 1). The total loads of vehicles 1,2,3 are respectively equal to 3 (serviced tasks 7,8,9), 2 (tasks 1,2) and 4 (tasks 3,4,5,6).

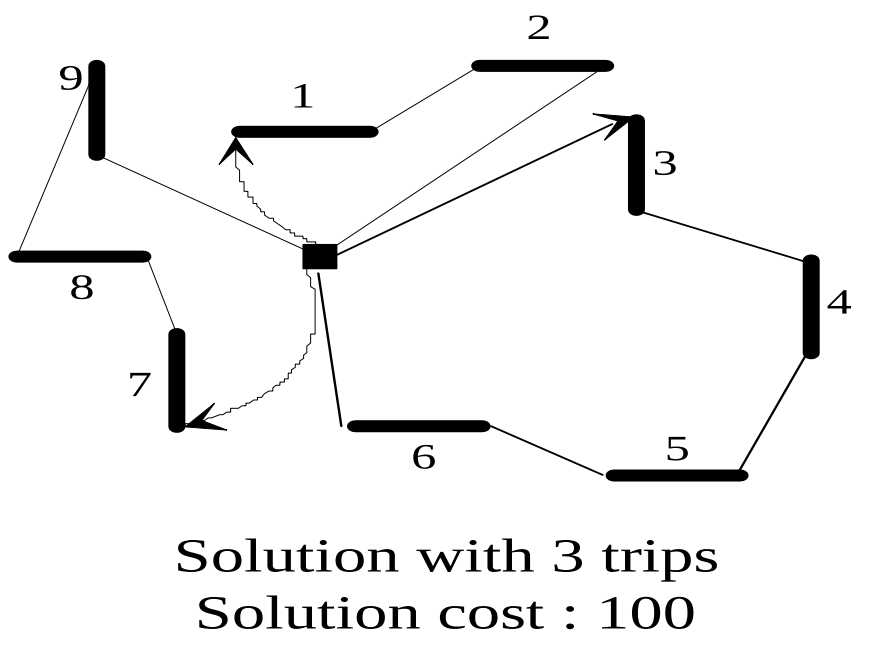

*Figure 1. A solution for the SCARP*

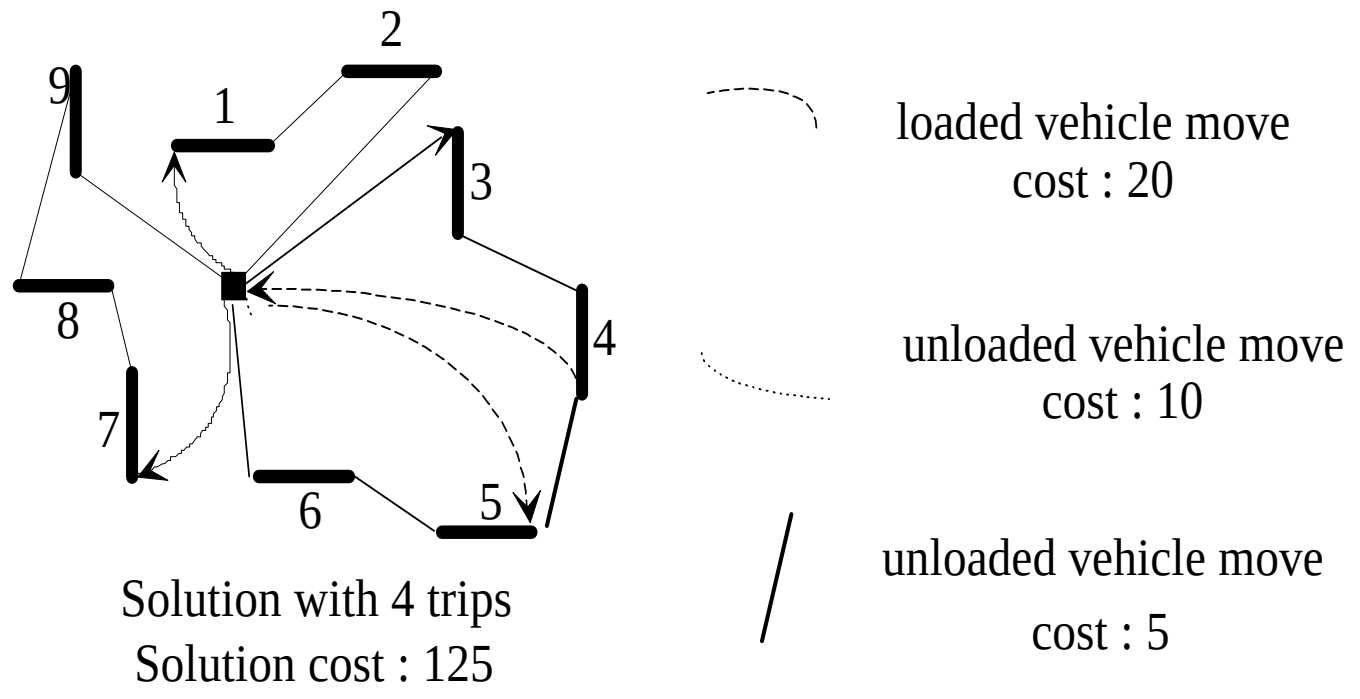

*Figure 2. The same solution for the SCARP*

Random events occurring in practice may affect estimated demands. Assume that after task 4, vehicle load is in fact 3.5 and the actual demands for tasks 5 and 6 become 1.2 and 1.4.

Because the vehicle can not service task 5, it moves from its current position (end of task 4) to the depot node and moves back to the beginning node of the task 5 to complete its trip (figure 2). This trip can be viewed as an additional trip. Due to this unproductive move, the solution cost becomes $H(S,\omega) = 125$ when $h(S) = 100$ . And one has $T(S,\omega) = 4$ when $t$(S)=3.

The objective of the SCARP consists in determining solutions close to the optimal solution of the CARP and robust in front of the random quantifies variations. (Fleury *et al.*, 2005) proposed a heuristic approach, (Fleury et al., 2004) proposed a memetic algorithm for determining robust solutions and (Lacomme *et al.*, 2005) proposed a memetic algorithm to minimize both average cost and standard deviation cost. The computational experiments prove that high quality solutions can be obtained in computational time not far from the computational time of the CARP.

### 1.3. The Multi-Objective CARP

The single objective CARP only deals with: minimizing the total cost of the trips. In fact, most waste management companies are also interested in balancing the trips. For instance, in Troyes (France), all 19 trucks leave the depot at 6 am and the waste collection must be completed as soon as possible to assign the crews to other tasks, e.g. sorting the waste at a recycling facility. Hence, the company wishes to solve a bi-objective version of the CARP, in which both the total duration of the trips and the duration of the longest trip (the makespan in scheduling theory) are to be minimized. This bi-objective CARP has been investigated first in (Lacomme *et al.*, 2005) whose proposed a non-dominated NSGA-II framework.

## 2. PROPOSAL FOR A STOCHASTIC MULTI-OBJECTIVE SCHEME FOR THE STOCHASTIC CARP

### 2.1. A framework for SCARP resolution

A framework for SCARP is here composed of two steps (figure 3): an optimization step and an evaluation of the robustness solutions (this second step only being used to evaluate the robustness of the best solution found at the end of the optimization process, then it can be omitted for the effective implementations): the best solution found is submitted to a replication phase consisting in statistic evaluations of robustness criteria in front of trials of random demands.

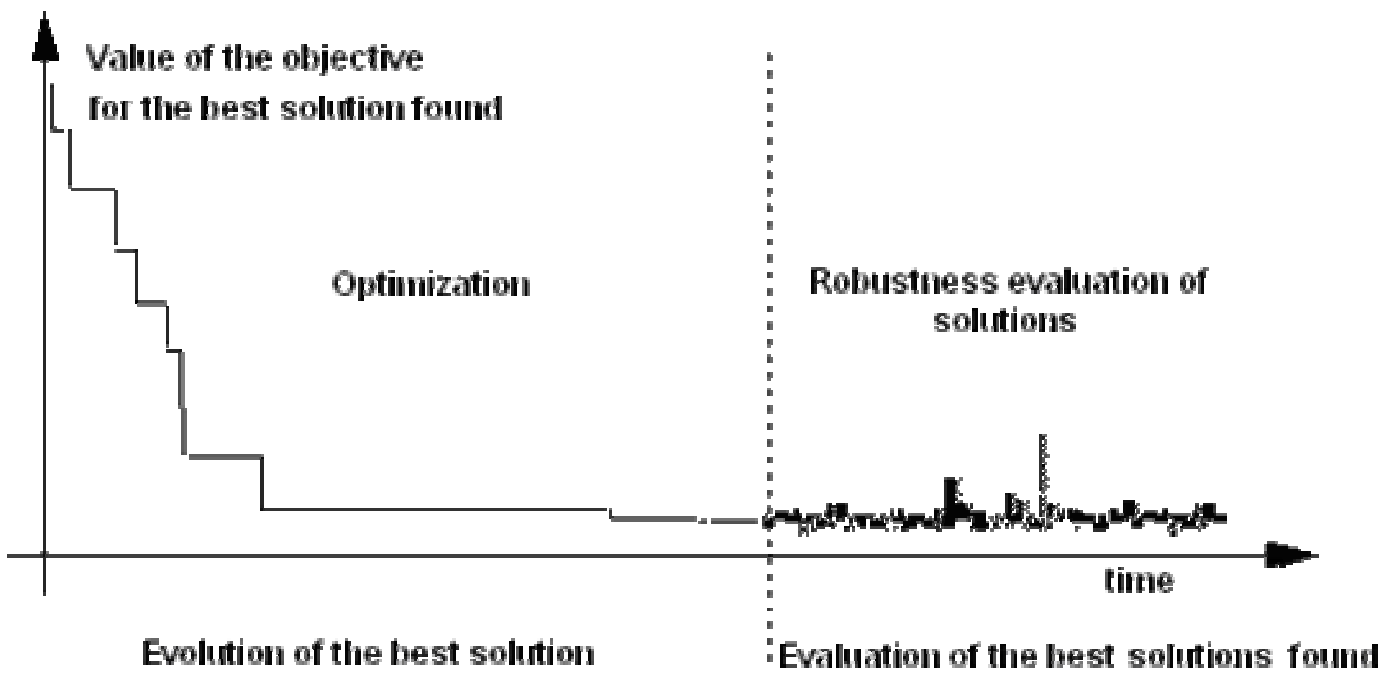

*Figure 3. Principle of evaluation of the stochastic problem resolution.*

### 2.2. First phase: the optimization phase

*Typical scheme of iterative methods for minimization of a stochastic criterion*

We address the wide-ranging problem of minimizing a random variable on a finite set. Let us precise the context:

- $X$ is a finite set,
- any $x$ of $X$ is associated to a random variable $H(x,\omega)$.

The problem to minimize $H$ on $X$ is not well defined, for " $x_1 < x_2 \; if \; H(x_1,\omega) < H(x_2,\omega)$ " is a random relation. The usual way is to transform this problem into a deterministic one, replacing the random variable $H$ by a deterministic one h. H usually depends on a lot of random variables $Q_i$ . Then h is evaluated with the same formula but where the random variables $Q_i$ are replaced, for example, by their expectation $q_i$ or by any deterministic value depending on the random variables $Q_i$ .

Thus the deterministic function $h$ is to be minimized. This technique may be irrelevant when non linear effects can strongly modify the value of $H$ . For example, when random events occur, a solution $x$ can become unrealistic *i.e.* $H(x,\omega)$ can be infinite (or very large) with a positive probability. When non linear effects exist, a better way is to replace (if possible) $H(x,\omega)$ itself by its expectation $\overline{H(x)}$, or another deterministic quantity associated to *H*. But when *H* can become infinite (or very high), it can be useful to search solutions *x* so that the probability $P\{H=\infty\}$ is at the most a given value ($P\{H=\infty\}\leq\varepsilon$). Other various criteria can of course be minimized. For example *h* can be either the conditional expectation of *H* given $\{H<\infty\}$ when $P\{H=\infty\}\leq\varepsilon$ and infinite when $P\{H=\infty\}>\varepsilon$. In these last cases, the first idea is to estimate *h* for example by $\overline{H(x,n)}$ (based on $n$ trials for the same solution $x$ with $n$ large enough, multiplying the minimization time nearly by $n$.)

The outline of such an iterative search process (Fleury, 1993) is described in figure 4. Nearly all stochastic metaheuristics which come from simulated annealing, taboo search… have extensions to tackle minimization of stochastic functions (Tsutsui and Ghosh, 1997) (Branke, 1998) (Ben-Tal and Nemirowski, 1998). For an introduction to stochastic scheduling and neighborhood based robustness approaches for scheduling, it is possible to refer to (Jensen, 2001).

1. Compute an initial solution $x_1$
2. Compute $\overline{H(x_1,n)}$
2. **Repeat**
3.1. Generate a solution $x_2$
3.2. Compute $\overline{H(x_2,n)}$
3.3. **If** $\overline{H(x_2,n)}\leq\overline{H(x_1,n)}$ **Then** $x_1 \leftarrow x_2$
3.4. **EndIf**
4. **Until** (Stop Condition)
5. **Return** $x_1$

*Figure 4. Outline of a basic iterative search process for a stochastic minimization*

However, proving the convergence of such a process is a challenging problem due to convergence conditions which highly depend on the function to minimize and on the generation of intermediate solutions. Fleury in 1993 (Fleury, 1993) has promoted an extension of the previous typical scheme in which the number of replications used to evaluate $\overline{H(x,n)}$ increases over the iterations of the algorithm. This extension permits to decrease the probability of error in accepting a new solution $x_2$ more promising than $x_1$. A demonstration in probability is proposed proving that under non restrictive hypotheses on the function to minimize the iterative process converges (with the probability one) towards robust solutions.

When an exact calculation or a high quality approximation of $\overline{H}$ can be mathematically performed, the minimization time is then reduced. Let us remark that mathematical analysis avoids errors due to the randomness of $\overline{H(x,n)}$. The line of research we here promote consists in coupling a mathematical analysis of solutions to a dedicated searching scheme devoted to the CARP.

*CARPs linked to the SCARP*

The only previous works which can be reported on the Stochastic CARP concern the tight and the slack approach of (Fleury *et al.*, 2004). The approach we promote here is quite different as regards both the objective and the models used:

- The *Law Approach*. This approach consists in the minimization of a deterministic objective depending on the laws of quantities to collect, for example $\overline{H(x)}$ or $\overline{H(x)}+k\sigma_H(x)$ (for a fixed *k*>0).

- The *tight approach* simply consists in solving the CARP linked to the SCARP using $\overline{Q_i}$ for the quantities to collect on the arc and the capacity $Q$ of vehicles. The function to minimize only depends on $\overline{Q_i}$, $Q$ and on the solution $x$. The experiments are fully available in (Fleury *et al.*, 2005).

 This approach can be denoted: $x \mapsto f(q_i, Q', x)$.

- The *slack approach* (Fleury *et al.*, 2004) is similar to the previous one but the optimization uses a smaller capacity of the vehicles $Q'$ ($Q' < Q$). The function to minimize is: $x \mapsto f(q_i, Q', x)$.

The *law approach* consists in solving the DCARP with a distribution law to modelize $Q_i$. The function to minimize ($x \mapsto g(\mathcal{L}(Q_i), Q, x)$) differs of *f* and here depends on the laws $\mathcal{L}(Q_i)$, on the vehicles capacity $Q$ and on the solution (table 1).

*Table 1. DCARP linked to the SCARP (here we denote the objective function by $f(Q_i, Q, S)$).*

| **SCARP – Mono-Objective resolution** | |
|---|---|
| • quantities to collect $Q_i$ are random variables<br>• capacity $Q$ of the vehicles is deterministic,<br>• $f(Q_i, Q, S)$ depends on the realization of the random variables $Q_i$ | |
| Associated **CARP** | **Comments** |
| Canonically associated CARP (*tight* approach) | $q_i = \overline{Q_i}$, $Q$ unchanged, objective $x \mapsto f(q_i, Q, x)$ |
| Associated CARP (*slack* approach) | $q_i = \overline{Q_i}$, $Q$ becomes $Q' < Q$, objective $x \mapsto f(q_i, Q', x)$ |
| Associated CARP (*law* approach) | $q_i = \overline{Q_i}$, $Q$ unchanged, objective function becomes $x \mapsto g(\mathcal{L}(Q_i), Q, x)$ where $\mathcal{L}(Q_i)$ is the law of $Q_i$. |

## 2.3. A MULTI-OBJECTIVE FRAMEWORK FOR THE SCARP

The graph *G* describing the problem is converted into an entirely directed internal graph *H*. The nodes are dropped and an arc index is used. Shortest path costs are pre-computed in a matrix *D*. For any pair of arcs *(u,v)*, *D(u,v)* is the traversal cost of a shortest path from *u* to *v*.

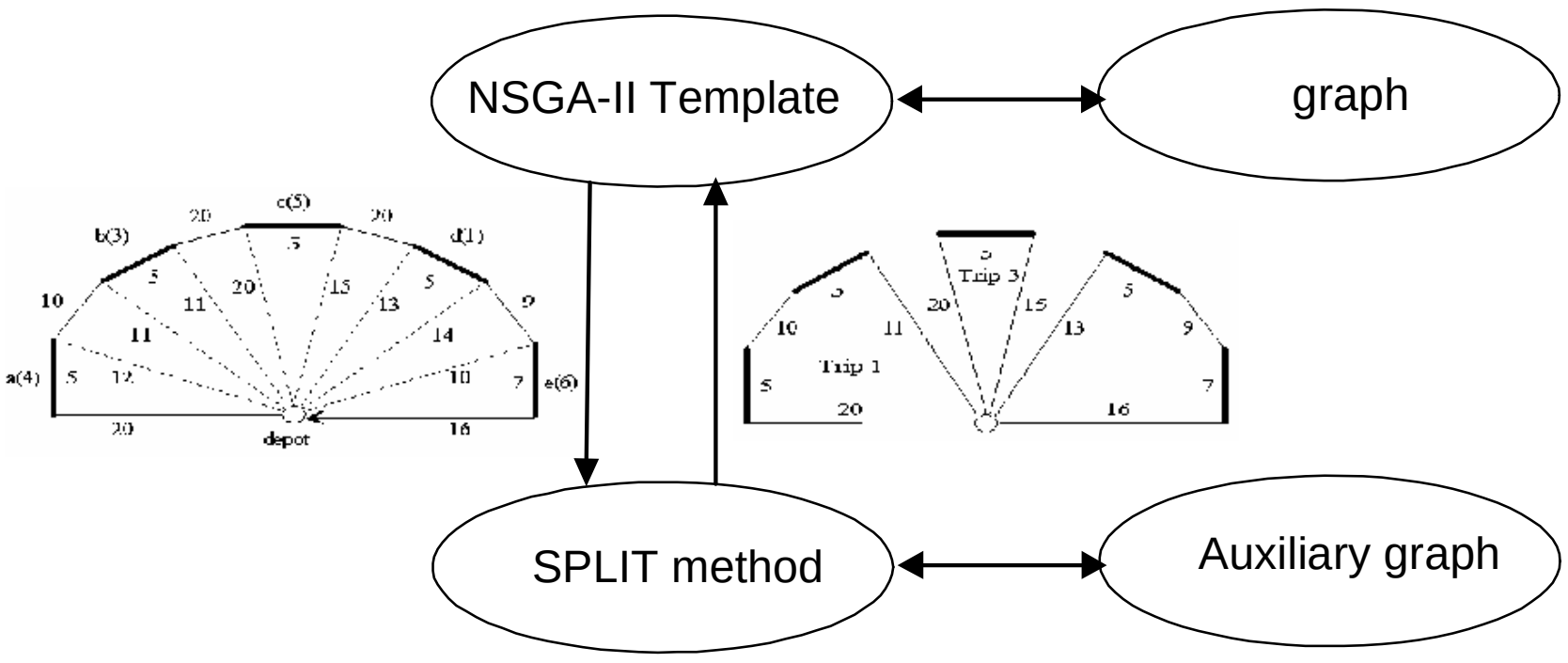

*Figure 4. NSGA-II template for SCARP resolution*

Today, several MOGA frameworks are available in literature and selecting the best one for a given problem is not obvious. A recent survey (Coello, 2000) and two comparative studies (Deb, 1999) (Zitzler *et al.*, 2000) try to provide guide-lines for selection, but these papers consider unconstrained problems, continuous objective functions, and specific sets of benchmarks. It is very difficult to draw conclusions for constrained combinatorial optimization problems.

(Lacomme *et al.*, 2005), finally turned their choice on the NSGA-II template to provide a multi-objective resolution of the CARP. A complete description of the NSGA-II template is available in (Lacomme *et al.*, 2005).

Each chromosome is an ordered set of required arcs assuming the same vehicle performs all trips in turn. This encoding is appealing because there always one optimal sequence. As stressed by (Lacomme *et al*, 2001) Ulusoy's algorithm provide a powerful technique for keeping a chromosome. The method consists in building a auxiliary graph in which each arc denotes a subsequence of required arcs. A shortest path algorithm in this graph gives the optimal split into trips taking into account the vehicles capacity (figure 5).

To tune this framework to the SCARP consist in defining the two stochastic criteria linked to the solution after the execution of the SPLIT method. The two criteria of interest are:

- the expected cost of a solution (expected total duration of the trips) and its standard deviation;
- the expected duration of the longest trip and its standard deviation.

Let us note $f_1(x)$ and $f_2(x)$ the two criteria to minimize.

## 2.4. A MULTI-OBJECTIVE FRAMEWORK FOR THE SCARP

The NSGA-II or Non-Dominated Sorting GA is an efficient multi-objective GA based on a non-dominated sorting of a population $P$ of $ns$ solutions. The sorting process starts by computing the non-dominated set of $P$ which defines the solutions of level 1 or front 1. Then, this set is temporarily removed from $P$ and the non-dominated set of the residual population is extracted to give the front 2 and son on until all solutions are classified as stressed in figure 5.

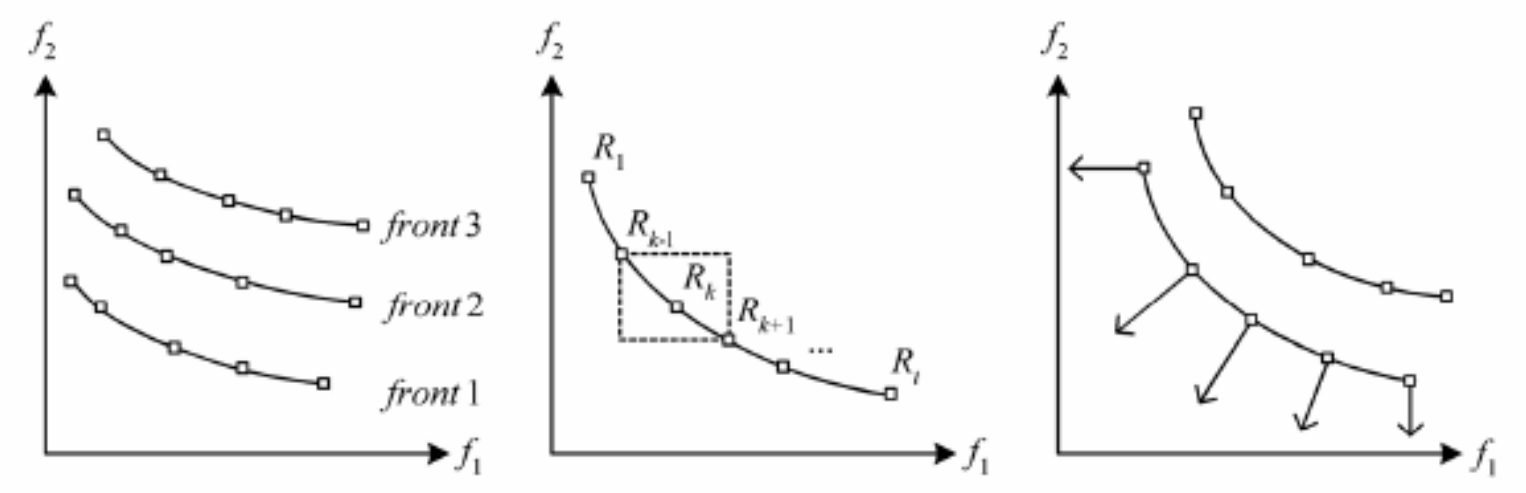

*Figure 5. NSGA-II template for SCARP resolution*

In NSGA-II, to each solution is assigned a fitness: its front number (1 being the best found value). To get well-spread fronts, parents are selected using a crowded tournament operator. A solution x wins a tournament with y if x has a better fitness (smaller front number) or if $x$ and $y$ are in the same front but $x$ has a larger crowding distance. This distance is depicted in the middle of figure 5 for a front $R$ of $t$ solutions. By convention, it is infinite for the two extreme points $R_1$ and $R_t$. For $R_k$, $1 < k < t$, it is equal to half of the perimeter of an enclosing rectangle with $R_{k-1}$ and $R_{k+1}$ placed on two vertices. This crowding distance is a kind of measure of the search space around $R_k$ which is not occupied by any other solution.

Starting from an initial population sorted by non-domination, one NSGA-II iteration consist of selecting $ns/2$ of parents with the crowded tournament operator, applying crossovers and mutations and adding the $ns$ resulting offspring's at the end of $P$, thus doubling its size. Finally, $P$ is reduced to its $ns$ best elements, using again a non-dominated sorting. An algorithm derived from NSGA-II is described in (Lacomme *et al.*, 2005) for the bi-objective CARP. Most components of the MA are recycled: the encoding of chromosomes, the OX-like crossover, the evaluation procedure Split and all the moves tested in the local search. The ways of integrating a local search (not foreseen in NSGA-II) without degrading solution dispersal are also studied. The best results are obtained by periodically applying to front 1 a local search with a direction depending on solution position, like on the right of figure 5.

More precisely, a move $x \to y$ is accepted if:

$$\pi_1.(f_1(y) - f_1(x)) + (1 - \pi_1).(f_2(y) - f_2(x)) < 0$$

with weight $\pi_1$ computed as in equation 1 where $f_k^{\min}$ (resp. $f_k^{\max}$ ) denotes the min (resp. max) value of the criterion $f_k$ , $k = 1,2$ .

$$\pi_1 = \left( \frac{f_1(x) - f_1^{\min}}{f_1^{\max} - f_1^{\min}} \right) \Bigg/ \left( \frac{f_1(x) - f_1^{\min}}{f_1^{\max} - f_1^{\min}} + \frac{f_2(x) - f_2^{\min}}{f_{12}^{\max} - f_2^{\min}} \right) \qquad (1)$$

The computational evaluation (Lacomme *et al.*, 2005) shows that, adding the local search sketched above provides better approximations of the Pareto set, while strongly accelerating convergence. Moreover, on a set of classical instances, the leftmost solution obtained at the end corresponds in most cases to the optimal or best-known solution of the single objective case.

### 2.5. MATHEMATICAL ANALYSIS OF SOLUTIONS

The challenging problem consists in applying the law approach simultaneously in two function denoted $f_1(x)$ and $f_2(x)$. These two deterministic functions must combine both an average value and standard deviation by positive penalties ( $\rho$ >0 and $\mu$ >0) in a linear combination:

$$f_2(x) = \overline{H(x)} + \rho.\sigma_H(x)$$

$$f_2(x) = \overline{M(x)} + \mu.\sigma_M(x)$$

where

$\overline{H(x)}$ is the expected cost of $x$ ;

$\sigma_H(x)$ is the standard deviation of the expected cost ;

$\overline{M(x)}$ is the expected duration of the longest trip of $x$ ;

$\sigma_M(x)$ is the standard deviation of the duration of the longest trip.

**Proposal for $\overline{H(x)}$ and $\sigma_H(x)$** (here we choose $q_i = \overline{Q_i}$ )

Considering the CARP associated to a SCARP, each trip $G_j$ of a solution $x = (G_i)_{1 \le i \le t(x)}$ satisfies: $\sum_{j \in G_i} q_j \le Q$ . Its (deterministic) cost is $\sum_{j \in G_i} c_j$ . Then the deterministic cost of the solution $x$ is $h(x) = \sum_{j=1}^{t(x)} \sum_{i \in G_j} c_i$ . In the SCARP, the cost $C_j$ of $G_j$ is greater than $\sum_{i \in G_j} c_i$ when $\sum_{i \in G_j} Q_i > Q$, because the vehicle must then go to the depot before continuing its task, thus $H(x,\omega) = \sum_{j=1}^{t(x)} C_j \ge h(x)$ . Hence, for any solution $x$ , $h(x)$ is the best possible value of $H(x,\omega)$ (let us recall that the cost of a trip only depends on the costs of arcs used but not on the collected quantities).

<u>Probability of additional move to the depot.</u>

For any trip $G_j$ of the CARP, the total amount of demand serviced is $\sum_{i \in G_j} q_i \le Q$. In the SCARP, as soon as $\sum_{i \in G_j} Q_i > Q$, the vehicle must turn back to the depot. The probability $p_j$ to report at most one move to the depot node during the trip $G_j$ of a solution $x$ is given by: $p_j = P\left\{\sum_{i \in G_j} Q_i > Q\right\}$. But $\sum_{i \in G_j} Q_i$ being a Gaussian random variable in the case of waste collection (thanks to the central limit theorem: the collected quantities along an arc is the sum of a large enough sum of the random quantities present in each container) with the expectation $\sum_{i \in G_j} q_i$ and the variance $\sigma_H\left(\sum_{i \in G_j} Q_i\right)^2 = \sum_{i \in G_j} \sigma(Q_i)^2$ (the random variables $Q_i$ are assumed to be independent, $\sigma(Q_i)$ being their standard deviation), so that $p_j = 1 - \varphi\left(\frac{Q - \sum_{i \in G_j} q_i}{\sqrt{\sum_{i \in G_j} \sigma(Q_i)^2}}\right)$ where $\varphi$ is the cumulative probability of *N(0,1)*: $\varphi(x) = \frac{1}{\sqrt{2.\pi}} \int_{-\infty}^{x} e^{-\frac{t^2}{2}}.dt$.

In the following, we assume (hypothesis 1) that the random variables $Q_i$ satisfy $\sigma(Q_i) = k.\overline{Q_i} = k.q_i$ for a fixed *k*>0, these formulas become respectively: $\sigma_H\left(\sum_{i \in G_j} Q_i\right)^2 = k^2.\sum_{i \in G_j} q_i^2$ and $p_j = 1 - \varphi\left(\frac{Q - \sum_{i \in G_j} q_i}{k.\sqrt{\sum_{i \in G_j} q_i^2}}\right)$. Hence, it is possible to evaluate, as soon as the demands are independent Gaussian random variables $N(q_{ik}, k^2 q_{ik}^2)$ the following characteristics:

- Probability of at most one additional trip in the solution is $P\{T(x,\omega) > t(x)\} = 1 - \prod_{j=1}^{t(x)} (1 - p_j)$ where $p_j$ is calculated as above.
- Probability of the most *m* additional trip in the solution is $P\{T(x,\omega) > t(x) + m\}$ (for a fixed integer $m \ge 1$) which is multinomial.

According to the capacity of the vehicle in front of the quantity to be collected on one arc, it seems reasonable to assume that at most one move can occur for any trip. This hypothesis (hypothesis 2) has been confirmed by the study (Fleury et al., 2005) and hence, the average number of trips is $\overline{T}(x) = t(x) + \sum_{j=1}^{t(x)} p_j$ and its standard deviation is $\sigma_T(x) = \sqrt{\sum_{j=1}^{t(x)} \left(p_j - p_j^2\right)}$.

<u>Position in a trip of an additional move to the depot.</u> Moreover, with a high probability, this additional move to the depot node will occur at the end of the trip, just before the last serviced arc of the trip (hypothesis 3). For a robust solution, the probability of an interruption is low, and therefore, it occurs with a high probability, just before the last serviced arc. In the following, we occasionally assume both following hypotheses are satisfied:

($H_2$) *any trip can be split into at the most two trips,*

($H_3$) *the additional move can occur only before the last serviced arc of the trip.*

Let us note $s_j$ the cost of an unproductive move from the last serviced arc of $G_j$ to the depot and from the depot to the next serviced arc of $G_j$. Under the hypotheses $(\mathrm{H}_1)$ and $(\mathrm{H}_2)$, with the probability $1- p_j$, the trip cost is $C_j = \sum_{i \in G_j} c_i$ and with the probability $p_j$ it is $C_j = \sum_{i \in G_j} (c_i + s_i)$. Under these hypotheses, it possible to compute some helpful characteristics of the random trip cost:

- the average cost of the trip: $\overline{H(G_j)} = \sum_{i \in G_j} (c_i + s_i . p_i)$.
- its standard deviation $\sigma_H(G_j) = \sqrt{s_j^2 .(c_i + s_i . p_i)}$.

For a solution $x = (G_j)_{1 \le j \le t(x)}$ composed of *t(x)* trips, assuming the demands are independent, the following properties can be established:

- The deterministic cost of solution is $h(x) = \sum_{j=1}^{t(x)} \sum_{i \in G_j} c_i$.
- The stochastic cost is $H(x, \omega) = \sum_{j=1}^{t(x)} C_j$.
- The average cost is $\overline{H(x)} = h(x) + \sum_{j=1}^{t(x)} s_j . p_j$.
- The standard deviation of the cost is $\sigma_H(x) = \sqrt{\sum_{j=1}^{t(x)} s_j^2 .(p_j - p_j^2)}$.

**Proposal for $\overline{M(x)}$ and $\sigma_M(x)$**

Let us, at first define the following deterministic numbers.

If no trip required a supplementary move to the depot node, $x_0 = \max_{1 \le i \le t(x)} C_i$.

Consider that only one trip requires a supplementary move to the depot node. For a trip $u$ $(1 \le u \le t(x))$,

$x_1^u = \max[C_u + s_u, x_0]$

f two trips $u$ and $v$ $(1 \le u,v \le t(x)$ and $u \ne v)$ require a supplementary move to the depot node,

$x_2^{uv} = \max[C_u + s_u, C_v + s_v, x_1^u]$

If three different trips $u$, $v$, $w$ $(1 \le u,v,w \le t)$ require a supplementary move to the depot node,

$x_3^{uvw} = \max[C_u + s_u, C_v + s_v, C_w + s_w, x_2^{uv}]$, and so on…

When all trips require a supplementary move to the depot, $x_{t(x)} = \max_{1 \le i \le t(x)} (C_i + s_i)$.

Thus, we have: $P\{X = x_0\} = \prod_{1 \le i \le t(x)} (1 - p_i) = \pi_0$ the probability that no trip requires any supplementary move to the depot:

$P\{X = x_1^u\} = p_u \cdot \prod_{\substack{1\le i\le t(x)\\ i\ne u}} (1-p_i) = \frac{p_u.\pi_0}{1-p_u}$ the probability that exactly one trip requires a supplementary move to the depot node;

$P\{X = x_2^{uv}\} = p_u \cdot p_v \cdot \prod_{\substack{1\le i\le t(x)\\ i\ne u\\ i\ne v}} (1-p_i) = \frac{p_u \cdot p_v.\pi_0}{(1-p_u)(1-p_v)}$ the probability that exactly two trips require a supplementary move to the depot node

$P\{X = x_3^{uvw}\} = p_u \cdot p_v \cdot p_w \cdot \prod_{\substack{1\le i\le t(x)\\ i\ne u\\ i\ne v\\ i\ne w}} (1-p_i) = \frac{p_u \cdot p_v \cdot p_w.\pi_0}{(1-p_u)(1-p_v)(1-p_w)}$ the probability that exactly three trips require a supplementary move to the depot node

And finally, $P\{X = x_{t(x)}\} = \prod_{1\le i\le t(x)} p_i$ is the probability that every trip requires a supplementary move to the depot.

The average length of the longest trip of $x$ can then be computed:

$$\overline{M(x)} = x_0.\pi_0 + \sum_{1\le u\le t(x)} x_1^u.\frac{p_u.\pi_0}{(1-p_u)} + \sum_{u=2}^{t(x)}\sum_{v=1}^{u-1} x_2^{uv}.\frac{p_u \cdot p_v.\pi_0}{(1-p_u)(1-p_v)}$$
$$+\sum_{u=3}^{t(x)}\sum_{v=2}^{u-1}\sum_{w=1}^{v-1} x_3^{uvw}.\frac{p_u \cdot p_v \cdot p_w.\pi_0}{(1-p_u)(1-p_v)(1-p_w)} + \cdots + x_{t(x)} \cdot \prod_{1\le j\le t(x)} p_j$$

It is possible to also obtain:

$$\overline{M(x^2)} = (x_0)^2.\pi_0 + \sum_{1\le u\le t(x)} (x_1^u)^2.\frac{p_u.\pi_0}{(1-p_u)} + \sum_{u=2}^{t(x)}\sum_{v=1}^{u-1} (x_2^{uv})^2.\frac{p_u \cdot p_v.\pi_0}{(1-p_u)(1-p_v)}$$
$$+\sum_{u=3}^{t(x)}\sum_{v=2}^{u-1}\sum_{w=1}^{v-1} (x_3^{uvw})^2.\frac{p_u \cdot p_v \cdot p_w.\pi_0}{(1-p_u)(1-p_v)(1-p_w)} + \cdots + (x_{t(x)})^2.\prod_{1\le j\le t(x)} p_j .$$

The standard deviation can then be computed: $\sigma_M(x) = \sqrt{\overline{M(x^2)} - \left(\overline{M(x)}\right)^2}$

<u>Implementation of $\overline{M(x)}$ and $\sigma_M(x)$</u>

Because $x_k^{u..z}$ is greater than $x_0$ and less than $x_{t(x)}$, one has

$$\overline{M(x)} \ge \pi_0.\left[x_0 + \sum_{u=2}^{t(x)} x_1^u.\frac{p_u}{(1-p_u)} + \sum_{u=2}^{t(x)}\sum_{v=1}^{u-1} x_2^{uv}\frac{p_u \cdot p_v}{(1-p_u)(1-p_v)}\right] + x_0.\left[1-\pi_0 - \sum_{u=2}^{t(x)}\frac{p_u\pi_0}{(1-p_u)} - \sum_{u=2}^{t(x)}\sum_{v=1}^{u-1}\frac{p_u \cdot p_v.\pi_0}{(1-p_u)(1-p_v)}\right]$$

and

$$\overline{M(x)} \le \pi_0.\left[x_0 + \sum_{u=2}^{t(x)} x_1^u.\frac{p_u}{(1-p_u)} + \sum_{u=2}^{t(x)}\sum_{v=1}^{u-1} x_2^{uv}\frac{p_u \cdot p_v}{(1-p_u)(1-p_v)}\right] + x_{t(x)}.\left[1-\pi_0 - \sum_{u=2}^{t(x)}\frac{p_u\pi_0}{(1-p_u)} - \sum_{u=2}^{t(x)}\sum_{v=1}^{u-1}\frac{p_u \cdot p_v.\pi_0}{(1-p_u)(1-p_v)}\right].$$

One can of course inprove these approximations using:

$$x_0.\sum_{u=3}^{t(x)}\sum_{v=2}^{u-1}\sum_{w=1}^{v-1}\frac{p_u \cdot p_v \cdot p_w.\pi_0}{(1-p_u)(1-p_v)(1-p_w)} \le \sum_{u=3}^{t(x)}\sum_{v=2}^{u-1}\sum_{w=1}^{v-1} x_3^{uvw}.\frac{p_u \cdot p_v \cdot p_w.\pi_0}{(1-p_u)(1-p_v)(1-p_w)} \le x_{t(x)}\sum_{u=3}^{t(x)}\sum_{v=2}^{u-1}\sum_{w=1}^{v-1}\frac{p_u \cdot p_v \cdot p_w.\pi_0}{(1-p_u)(1-p_v)(1-p_w)}$$

and so on…

Similar considerations are valid for $\overline{M(x^2)}$.

The previous remarks show that, if $e \le \overline{M(x)} \le E$ and $c \le \overline{M(x^2)} \le C$, then $\sqrt{\max\left(0, c - E^2\right)} \le \sigma_M(x) \le \sqrt{C - e^2}$ and then:

$\overline{M(x)} \approx \frac{e+E}{2}$ and $\sigma_M(x) \approx \frac{\sqrt{\max\left(0, c-E^2\right)} + \sqrt{C - e^2}}{2}$.

## 2.6. EVALUATION OF ROBUSTNESS AND EVALUATION OF SOLUTIONS QUALITY

**Robustness of solutions**

The second phase consists in gathering statistics. Let $x$ be one of the best solutions obtained at the end of the optimization phase. $n$ replications can be performed for a careful analysis of solution properties as regards robustness criteria. The statistics can include (but are not limited to):

- $\overline{H(x,n)}$ : the average cost over $n$ independent evaluations of $H(x,\omega)$. This estimates the expectation $\overline{H(x)}$.
- $\overline{M(x,n)}$ : the average cost over $n$ independent evaluations of $M(x,\omega)$. This estimates the expectation $\overline{M(x)}$
- $\sqrt{\frac{n}{n-1}}.\sigma_H(x,n)$ where $\sigma_H(x,n)$ is the standard deviation of the cost over $n$ independent evaluations of $x$. This estimates the standard deviation $\sigma_H(x)$.
- $\sqrt{\frac{n}{n-1}}.\sigma_M(x,n)$ where $\sigma_M(x,n)$ is the standard deviation of the number of trips over $n$ independent evaluations of $x$. This estimates the standard deviation $\sigma_M(x)$.
- $E_{\overline{H(x)}} = \frac{\overline{H(x)} - \overline{H(x,n)}}{\overline{H(x,n)}} \times 100$ represents the gap (in percent) between $\overline{H(x,n)}$ and $\overline{H(x)}$ and $E_{\overline{M(x)}} = \frac{\overline{M(x)} - \overline{M(x,n)}}{\overline{M(x,n)}} \times 100$ represents the gap (in percent) between $\overline{M(x,n)}$ and $\overline{M(x)}$.
- $E_{\sigma_H(x)} = \frac{\overline{\sigma_H(x)} - \overline{\sigma_H(x,n)}}{\overline{H(x)}} \times 100$ is the relative importance (in percent) of the error in the standard deviation evaluation for $\overline{H(x)}$ and $E_{\sigma_M(x)} = \frac{\overline{\sigma_M(x)} - \overline{\sigma_M(x,n)}}{\overline{M(x)}} \times 100$ is the relative importance (in percent) of the error in the standard deviation evaluation for $\overline{M(x)}$

For convenience, we note $E^1_{\overline{H(x)}}$, $E^1_{\overline{M(x)}}$, $E^1_{\sigma_H(x)}$, $E^1_{\sigma_M(x)}$ the values of the leftmost solution and $E^2_{\overline{H(x)}}$, $E^2_{\overline{M(x)}}$, $E^2_{\sigma_H(x)}$, $E^2_{\sigma_M(x)}$ the values of the rightmost solution.

**Representation of solutions**

Representing the solutions, we propose (figure 6) to highlight the standard deviation of both objective functions (figure 6). The square represent solutions for which both the expected cost of a solution and the expected duration of the longest trip are included in the average value plus/minus the standard deviation. The size of the squares gives a graphical representation of solutions robustness: large size squares denote high sensitive solutions and small size squares denotes robust solutions. Mathematical expressions and replications evaluation can conduct to different evaluation of both criteria and of standard deviations.

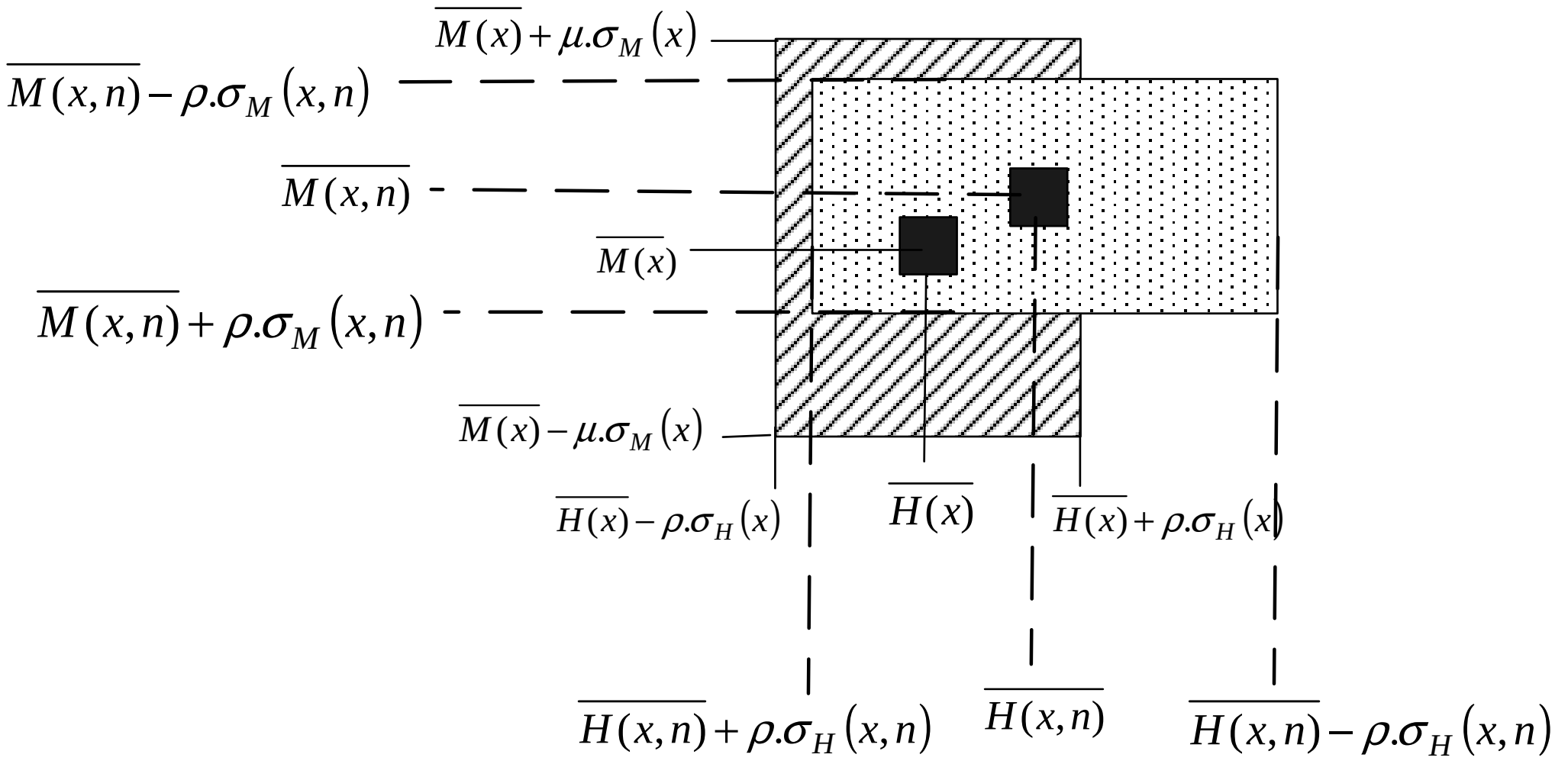

*Figure 6. Representation of one solution as square depending on the"exact" or evaluated expressions.*

**Quality of evaluations of the solutions**

Let $\left(\overline{H(x)}^1;\overline{M(x)}^1\right)$ be the leftmost solution of the first front and $\left(\overline{H(x)}^2;\overline{M(x)}^2\right)$ be the rightmost solution of the first front. These two solutions can be respectively compared to the leftmost solution $\left(h(x)^1;m(x)^1\right)$ and the rightmost one $\left(h(x)^2;m(x)^2\right)$ (Lacomme *et al.*, 2005) obtained by solving the bi-objective (deterministic) CARP. Both solutions can also be compared to the best solution $\overline{H(x)}$ found solving the (stochastic mono-objective) SCARP (Fleury *et al.*, 2004) (Fleury et al., 2005). To avoid any ambiguity, this value is denoted $\overline{H(x)}^{Mono}$ in the rest of the paper. Of course, neither $\left(h(x)^1;m(x)^1\right)$ nor $\left(h(x)^2;m(x)^2\right)$ are a lower bound, but the best found solution solving the CARP using a bi-objective optimization scheme. Figure 7 gives a graphical representation of all the solutions available depending on the resolution scheme applied.

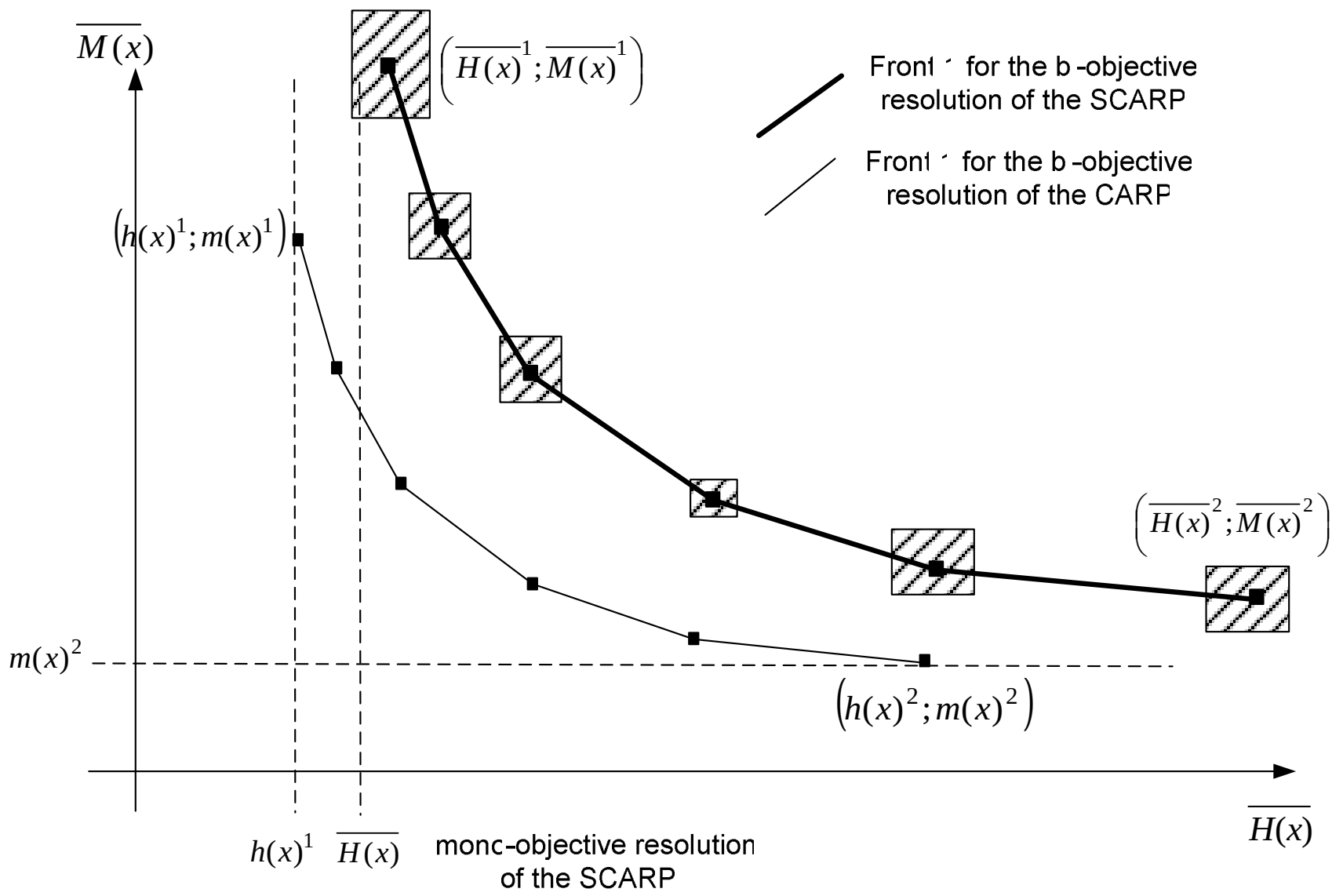

*Figure 7. Evaluation of solutions for the bi-objective resolution of the SCARP as regard previous published results on the SCARP (mono-objective) and on the CARP (bi-objective resolution).*

The leftmost and rightmost solutions obtained using the bi-objective resolution scheme of the SCARP, can be compared to the leftmost and the rightmost solutions found using the bi-objective resolution of the (deterministic) CARP. Let us note $E^1_{\overline{H(x)}^1-h(x)^1}$, $E^1_{\overline{M(x)}^1-m(x)^1}$, $E^2_{\overline{H(x)}^2-h(x)^2}$ and $E^2_{\overline{M(x)}^2-m(x)^2}$ the gap (in percent) from rightmost and leftmost solutions. These gaps can be computed according to the next formula:

$$E^1_{\overline{H(x)}^1-h(x)^1}=\left(\left(\overline{H(x)}^1-h(x)^1\right)\Big/h(x)^1\right)\times 100 \qquad E^2_{\overline{H(x)}^2-h(x)^2}=\left(\left(\overline{H(x)}^2-h(x)^2\right)\Big/h(x)^2\right)\times 100$$

$$E^1_{\overline{M(x)}^1-m(x)^1}=\left(\left(\overline{M(x)}^1-m(x)^1\right)\Big/m(x)^1\right)\times 100 \qquad E^2_{\overline{M(x)}^2-m(x)^2}=\left(\left(\overline{M(x)}^2-m(x)^2\right)\Big/m(x)^2\right)\times 100$$

The leftmost solution of the SCARP resolution by the bi-objective scheme can be compared with the solution obtained using a mono-objective resolution scheme of the SCARP. This gap can be computed thinks to:

$$E^1_{\overline{H(x)}^1-\overline{H(x)}^{MONO}}=\left(\left(\overline{H(x)}^1-\overline{H(x)}^{MONO}\right)\Big/\overline{H(x)}^{MONO}\right)\times 100$$

Figure 8 provides a graphical representation of solutions quality criteria as regards previous published results on the mono-objective resolution of the SCARP and on the bi-objective resolution of the CARP.

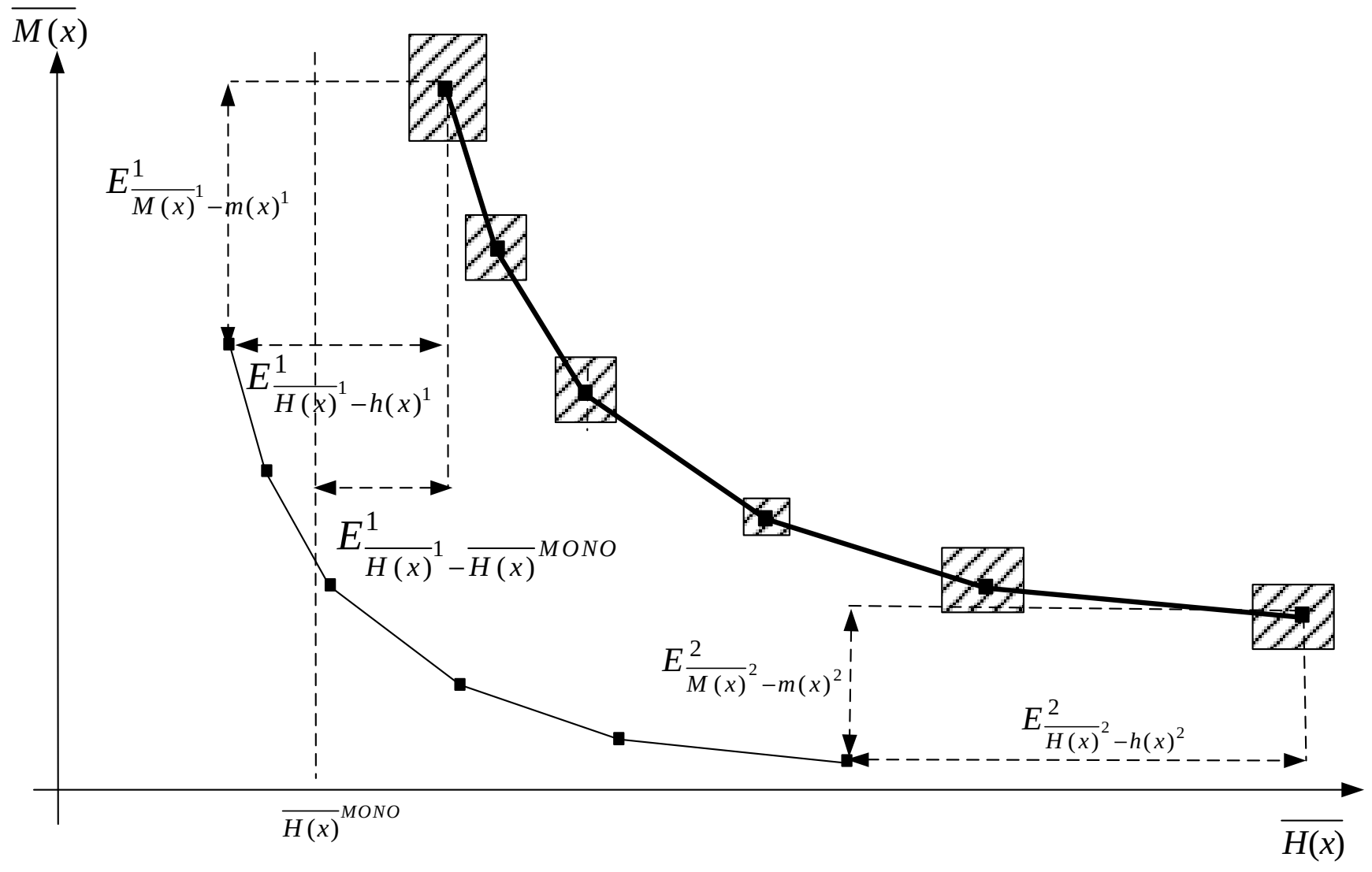

*Figure 8. Graphical representation of solutions quality criteria*

## 3. NUMERICAL EXPERIMENTS

The experiments have been carried out using the well-known instances introduced by (Belenguer and Benavent, 2003), (Eglese and Li, 1996) and (DeArmon, 1981). These computational evaluations have been achieved on a Pentium IV 2.8 GHz with 512Mo using Windows XP operating system. The program has been developed using Delphi 7.0 package. The demands $q_i$ are replaced by independent truncated random variables $N\left(q_i, r^2.q_i^2\right)$ with $r = 0.1$. The penalties $\rho$ and $\mu$ in the linear combination used for $f_1$ and $f_2$ are equal to 10. All the experiments have been carried out with a population of 60 chromosomes, 1000 iterations and a directed local search applied every 10 iterations. This set of parameters has been applied for all the experiments. When the first front was composed of one solution only, this solution is considered for being both leftmost and rightmost solution.

### 3.1. ROBUSTNESS OF SOLUTIONS FOR THE GDB INSTANCES

The solutions can be represented as points in a set of fronts (figure 9) or can be represented as squares in a set of fronts (figure 10) for the initial population.

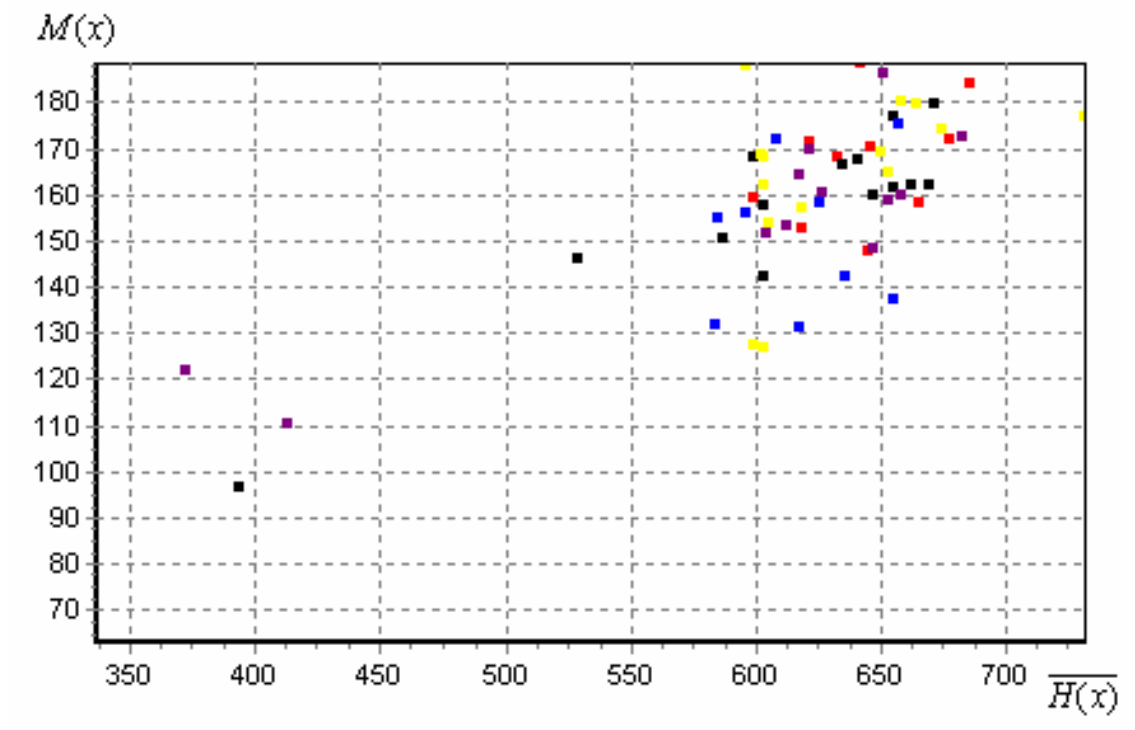

*Figure 9. Representation of the initial solutions of Gdb1 instance*

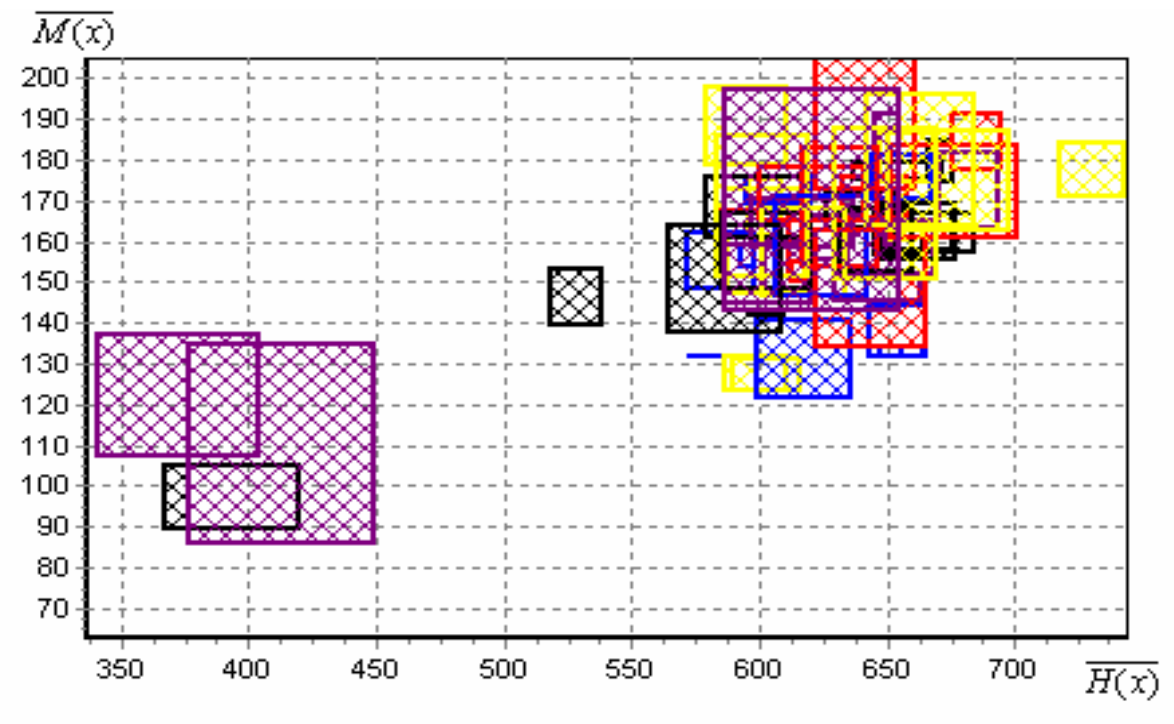

*Figure 10. Representation by domains of the initial solutions of Gdb1 instance*

One can note, that the second representation permits to show if the considered solutions have or not a great standard deviation (figure 11).

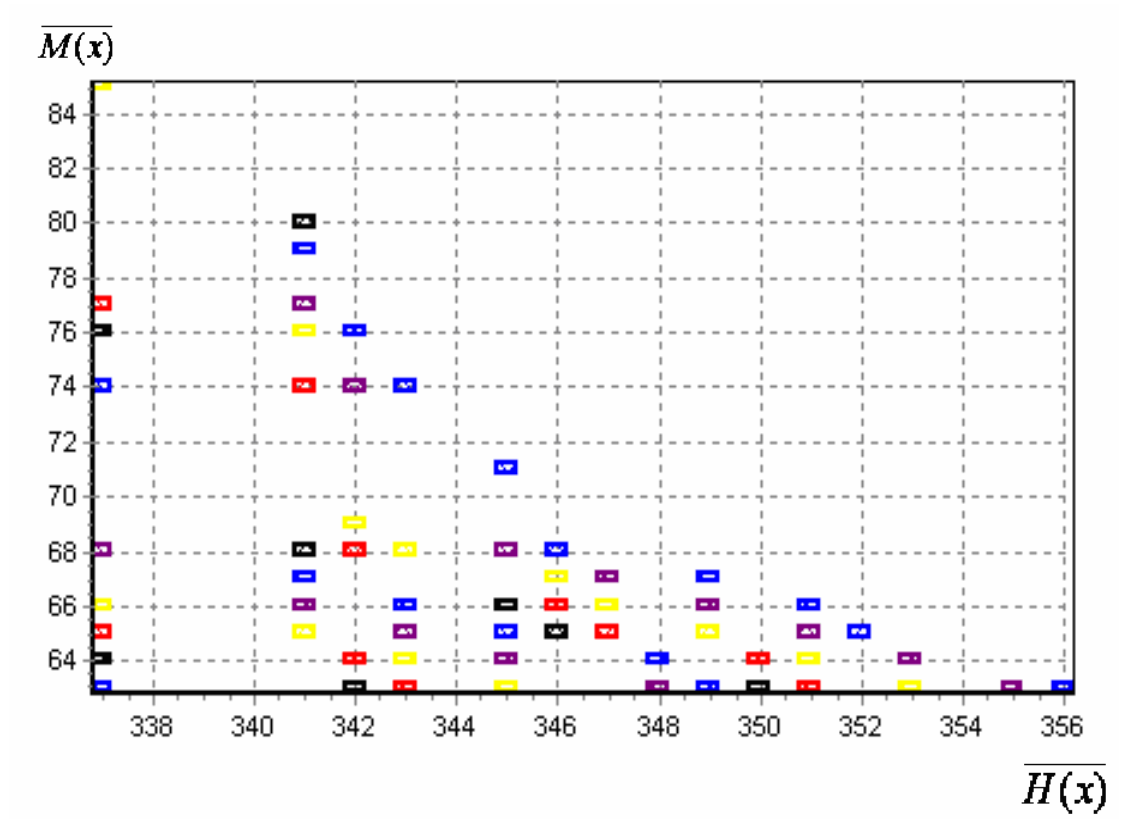

*Figure 11. Representation by domains of the solutions of Gdb1 instance after 1000 iterations*

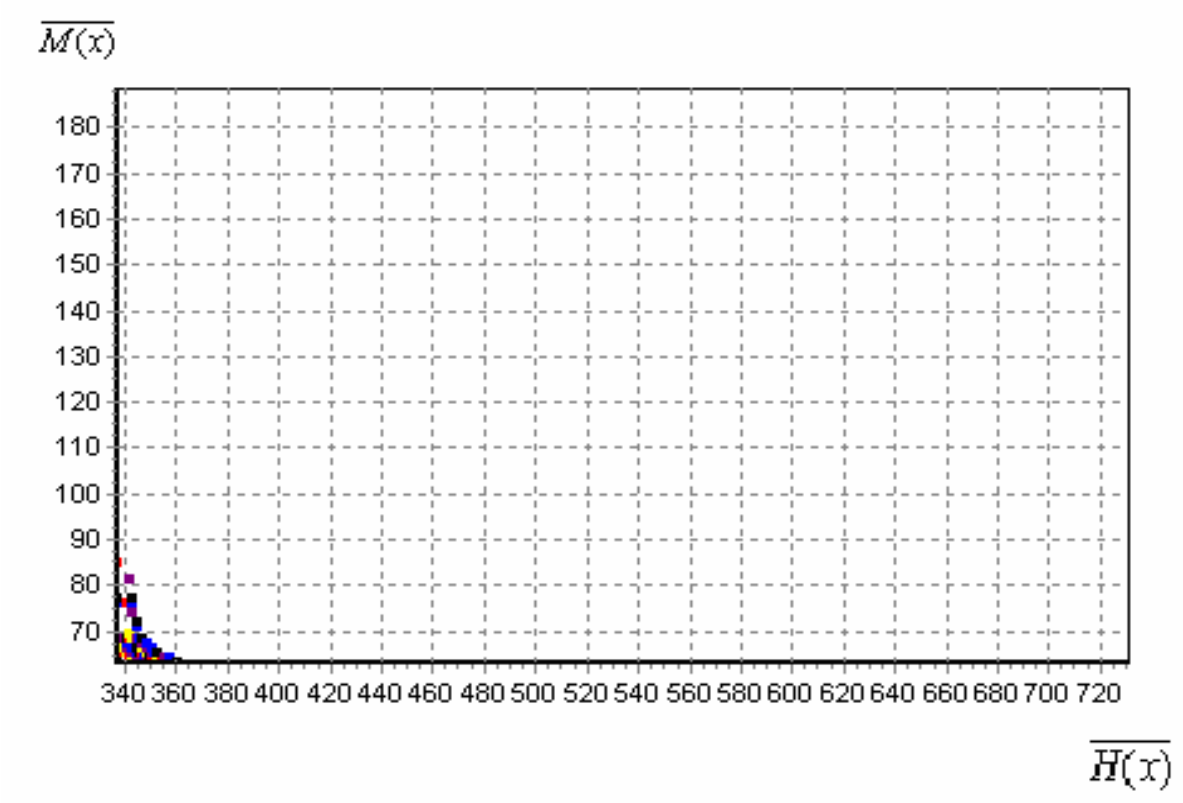

*Figure 12. Representation by domains of the solutions of Gdb1 instance after 1000 iterations – Same scale than figure 9*

Table 2 presents an analysis of solutions for Gdb instances. The aim is to highlight the high correlation between the mathematical evaluation of solutions and the evaluation provides by replications. The evaluation concerns the rightmost and the leftmost solutions obtained at the end of the optimization phase.

*Table 2. Mathematical evaluation and replications evaluation of the rightmost and of the leftmost solutions of Gdb instances*

| | **Numerical values provided by mathematical expressions** | | | | **Numerical values provided by replications phase** | | | |
|---|---|---|---|---|---|---|---|---|
| **Gdb** | $\overline{H(x)}$ | $\overline{M(x)}$ | $\sigma_H(x)$ | $\sigma_M(x)$ | $\overline{H(x,n)}$ | $\overline{M(x,n)}$ | $\sigma_H(x,n)$ | $\sigma_M(x,n)$ |
| **1** | 337.0 | 63.0 | 0.0 | 0.0 | 337.0 | 63.0 | 0.0 | 0.0 |
| | 337.0 | 63.0 | 0.0 | 0.0 | 337.0 | 63.0 | 0.0 | 0.0 |
| **2** | 367.0 | 62.0 | 0.0 | 0.0 | 367.0 | 62.0 | 0.0 | 0.0 |
| | 417.0 | 59.0 | 0.0 | 0.0 | 417.0 | 59.0 | 0.0 | 0.0 |
| **3** | 296.0 | 60.0 | 0.0 | 0.0 | 296.0 | 60.0 | 0.0 | 0.0 |
| | 347.0 | 59.0 | 0.0 | 0.0 | 347.0 | 59.0 | 0.0 | 0.0 |
| **4** | 317.0 | 72.0 | 0.0 | 0.0 | 317.0 | 72.0 | 0.0 | 0.0 |
| | 350.0 | 64.0 | 0.0 | 0.0 | 350.0 | 64.0 | 0.0 | 0.0 |
| **5** | 433.0 | 68.0 | 0.0 | 0.0 | 433.0 | 68.0 | 0.0 | 0.0 |
| | 479.0 | 64.0 | 0.0 | 0.0 | 479.0 | 64.0 | 0.0 | 0.0 |
| **6** | 324.0 | 68.0 | 0.0 | 0.0 | 324.0 | 68.0 | 0.0 | 0.0 |
| | 351.0 | 64.0 | 0.0 | 0.0 | 351.0 | 64.0 | 0.0 | 0.0 |
| **7** | 359.0 | 68.0 | 0.0 | 0.0 | 359.0 | 68.0 | 0.0 | 0.0 |
| | 368.0 | 66.0 | 0.0 | 0.0 | 368.0 | 66.0 | 0.0 | 0.0 |
| **8** | 399.1 | 44.0 | 0.6 | 0.1 | 399.1 | 44.0 | 0.6 | 0.0 |
| | 493.0 | 39.0 | 0.0 | 0.0 | 493.0 | 39.0 | 0.0 | 0.0 |
| **9** | 361.2 | 38.0 | 0.9 | 0.1 | 361.2 | 38.0 | 0.7 | 0.0 |
| | 379.1 | 37.0 | 0.7 | 0.0 | 379.1 | 37.0 | 0.6 | 0.0 |
| **10** | 283.0 | 66.0 | 0.0 | 0.0 | 283.0 | 66.0 | 0.0 | 0.0 |
| | 329.0 | 57.0 | 0.0 | 0.0 | 329.0 | 57.0 | 0.0 | 0.0 |
| **11** | 409.0 | 81.0 | 0.0 | 0.0 | 409.0 | 81.0 | 0.0 | 0.0 |
| | 439.0 | 61.0 | 0.0 | 0.0 | 439.0 | 61.0 | 0.0 | 0.0 |
| **12** | 523.7 | 99.1 | 3.9 | 2.6 | 523.6 | 99.1 | 3.9 | 2.7 |
| | 613.1 | 93.0 | 1.6 | 0.0 | 613.1 | 93.0 | 1.7 | 0.0 |
| **13** | 556.0 | 128.0 | 0.1 | 0.0 | 556.0 | 128.0 | 0.0 | 0.0 |
| | 556.0 | 128.0 | 0.1 | 0.0 | 556.0 | 128.0 | 0.0 | 0.0 |
| **14** | 104.0 | 20.0 | 0.0 | 0.0 | 104.0 | 20.0 | 0.0 | 0.0 |
| | 108.0 | 17.0 | 0.0 | 0.0 | 108.0 | 17.0 | 0.0 | 0.0 |
| **15** | 58.0 | 16.0 | 0.0 | 0.0 | 58.0 | 16.0 | 0.0 | 0.0 |
| | 60.0 | 13.0 | 0.0 | 0.0 | 60.0 | 13.0 | 0.0 | 0.0 |
| **16** | 131.0 | 20.0 | 0.0 | 0.0 | 131.0 | 20.0 | 0.0 | 0.0 |
| | 135.0 | 18.0 | 0.0 | 0.0 | 135.0 | 18.0 | 0.0 | 0.0 |
| **17** | 91.0 | 17.0 | 0.0 | 0.0 | 91.0 | 17.0 | 0.0 | 0.0 |
| | 91.0 | 17.0 | 0.0 | 0.0 | 91.0 | 17.0 | 0.0 | 0.0 |
| **18** | 172.0 | 33.0 | 0.0 | 0.0 | 172.0 | 33.0 | 0.0 | 0.0 |
| | 182.0 | 27.0 | 0.0 | 0.0 | 182.0 | 27.0 | 0.0 | 0.0 |
| **19** | 63.0 | 17.0 | 0.0 | 0.0 | 63.0 | 17.0 | 0.0 | 0.0 |
| | 63.0 | 18.0 | 0.0 | 0.0 | 63.0 | 18.0 | 0.0 | 0.0 |
| **20** | 123.0 | 33.0 | 0.0 | 0.0 | 123.0 | 33.0 | 0.0 | 0.0 |
| | 127.0 | 21.0 | 0.0 | 0.0 | 127.0 | 21.0 | 0.0 | 0.0 |
| **21** | 160.0 | 28.0 | 0.0 | 0.0 | 160.0 | 28.0 | 0.0 | 0.0 |
| | 166.0 | 22.0 | 0.0 | 0.0 | 166.0 | 22.0 | 0.0 | 0.0 |
| **22** | 207.0 | 26.0 | 0.1 | 0.0 | 207.0 | 26.0 | 0.0 | 0.0 |
| | 208.0 | 20.0 | 0.0 | 0.0 | 208.0 | 20.0 | 0.0 | 0.0 |
| **23** | 239.0 | 28.0 | 0.1 | 0.0 | 239.0 | 28.0 | 0.0 | 0.0 |
| | 245.0 | 19.0 | 0.1 | 0.0 | 245.0 | 19.0 | 0.0 | 0.0 |

*Table 3. Error (in percentage) between mathematical expressions and replications evaluations for the rightmost and of the leftmost solutions*

| Gdb | $E_{\overline{H(x)}}$ | $E_{\overline{M(x)}}$ | $E_{\sigma_H(x)}$ | $E_{\sigma_M(x)}$ |
|---|---|---|---|---|
| **1** | 0.00 | 0.00 | 0.00 | 0.00 |
| | 0.00 | 0.00 | 0.00 | 0.00 |
| **2** | 0.00 | 0.00 | 0.00 | 0.00 |
| | 0.00 | 0.00 | 0.00 | 0.00 |
| **3** | 0.00 | 0.00 | 0.00 | 0.00 |
| | 0.00 | 0.00 | 0.00 | 0.00 |
| **4** | 0.00 | 0.00 | 0.00 | 0.00 |
| | 0.00 | 0.00 | 0.00 | 0.00 |
| **5** | 0.00 | 0.00 | 0.00 | 0.00 |
| | 0.00 | 0.00 | 0.00 | 0.00 |
| **6** | 0.00 | 0.00 | 0.00 | 0.00 |
| | 0.00 | 0.00 | 0.00 | 0.00 |
| **7** | 0.00 | 0.00 | 0.00 | 0.00 |
| | 0.00 | 0.00 | 0.00 | 0.00 |
| **8** | 0.00 | 0.00 | 0.00 | 0.14 |
| | 0.00 | 0.00 | 0.01 | 0.00 |
| **9** | 0.01 | 0.00 | 0.05 | 0.34 |
| | 0.01 | 0.00 | 0.03 | 0.00 |
| **10** | 0.00 | 0.00 | 0.00 | 0.00 |
| | 0.00 | 0.00 | 0.00 | 0.00 |
| **11** | 0.00 | 0.00 | 0.00 | 0.00 |
| | 0.00 | 0.00 | 0.00 | 0.00 |
| **12** | 0.02 | -0.01 | 0.00 | -0.18 |
| | 0.00 | 0.00 | -0.02 | 0.00 |
| **13** | 0.00 | 0.00 | 0.02 | 0.00 |
| | 0.00 | 0.00 | 0.00 | 0.00 |
| **14** | 0.00 | 0.00 | 0.00 | 0.00 |
| | 0.00 | 0.00 | 0.00 | 0.00 |
| **15** | 0.00 | 0.00 | 0.00 | 0.00 |
| | 0.00 | 0.00 | 0.00 | 0.00 |
| **16** | 0.00 | 0.00 | 0.00 | 0.00 |
| | 0.00 | 0.00 | 0.00 | 0.00 |
| **17** | 0.00 | 0.00 | 0.00 | 0.00 |
| | 0.00 | 0.00 | 0.00 | 0.00 |
| **18** | 0.00 | 0.00 | 0.00 | 0.00 |
| | 0.00 | 0.00 | 0.00 | 0.00 |
| **19** | 0.00 | 0.00 | 0.01 | 0.00 |
| | 0.00 | 0.00 | 0.01 | 0.00 |
| **20** | 0.00 | 0.00 | 0.00 | 0.00 |
| | 0.00 | 0.00 | 0.00 | 0.00 |
| **21** | 0.00 | 0.00 | 0.05 | 0.00 |
| | 0.00 | 0.00 | 0.00 | 0.00 |
| **22** | 0.00 | 0.00 | 0.03 | 0.00 |
| | 0.00 | 0.00 | 0.02 | 0.00 |
| **23** | 0.00 | 0.00 | 0.00 | 0.00 |
| | 0.00 | 0.00 | 0.00 | 0.00 |
| **Avg.** | **0.00** | **0.00** | **0.0** | **0.01** |
| Max | 0.02 | 0.00 | 0.05 | 0.34 |
| Min | 0.00 | -0.01 | -0.02 | 0.18 |

The equation between mathematical expressions and replications evaluations is In the table 4, for the 23 instances, over 46 solutions, the average

- $E_{\overline{H(x)}}$ is the gap between $\overline{H(x,n)}$ and $\overline{H(x)}$
- $E_{\overline{M(x)}}$ is the gap between $\overline{M(x,n)}$ and $\overline{M(x)}$.
- $E_{\sigma_H(x)}$ is the relative importance (in percent) of the error in the standard deviation evaluation for $\overline{H(x)}$
- $E_{\sigma_M(x)}$ is the relative importance (in percent) of the error in the standard deviation evaluation for $\overline{M(x)}$

The average gap between $\overline{H(x,n)}$ and $\overline{H(x)}$ remains 0 over the 46 solutions investigated. The same is true for the gap between $\overline{M(x,n)}$ and $\overline{M(x)}$, and for the relative importance of the error in standard deviation of $\overline{H(x)}$ and $\overline{M(x)}$.

Similar high quality results are also reported for the Val instances and for the Eglese instances. This study implies that the mathematical formulas provide high quality estimations of the four studied criteria.

In the next section, the quality of the solutions is compared to the previous published on the SCARP with a mono-objective approach and on the CARP with a bi-objective one (figure 8).

## 3.2. REMARKS ON THE MATHEMATICAL FORMULAE (VAL INSTANCES)

To obtain a suitable evaluation of mathematical formula performance, table 4a and table 4b give the error between mathematical formula and replications evaluation. One can note, that whatever the criteria and whatever the solution (leftmost or rightmost) the error is less than 0.02%. In appendix; table A1, gives all the values obtained for the standard deviation.

*Table 4a and Table 4b. Error between mathematical expressions and replications evaluations for the rightmost and of the leftmost solutions*

| | Leftmost solution | | Rightmost solution | |
|---|---|---|---|---|
| **Val** | $E^1_{\overline{H(x)}}$ | $E^1_{\overline{M(x)}}$ | $E^2_{\overline{H(x)}}$ | $E^2_{\overline{M(x)}}$ |
| **1a** | 0.00 | 0.00 | 0.00 | 0.00 |
| **1b** | 0.00 | 0.00 | 0.00 | 0.00 |
| **1c** | 0.00 | 0.00 | 0.00 | 0.00 |
| **2a** | 0.00 | 0.00 | 0.00 | 0.00 |
| **2b** | 0.00 | 0.00 | 0.00 | 0.00 |
| **2c** | -0.01 | 0.00 | 0.00 | 0.00 |
| **3a** | 0.00 | 0.00 | 0.00 | 0.00 |
| **3b** | -0.04 | -0.06 | 0.00 | 0.00 |
| **3c** | 0.00 | 0.00 | 0.00 | 0.00 |
| **4a** | 0.00 | 0.00 | 0.00 | 0.00 |
| **4b** | 0.00 | 0.00 | 0.00 | 0.00 |
| **4c** | 0.00 | 0.00 | 0.00 | 0.00 |
| **4d** | -0.15 | -1.08 | 0.00 | 0.00 |
| **5a** | 0.00 | 0.00 | 0.00 | 0.00 |
| **5b** | 0.00 | 0.00 | 0.00 | 0.00 |
| **5c** | 0.00 | 0.00 | 0.00 | 0.00 |
| **5d** | 0.00 | 0.00 | 0.00 | 0.00 |
| **6a** | 0.00 | 0.00 | 0.00 | 0.00 |
| **6b** | 0.00 | 0.00 | 0.00 | 0.00 |
| **6c** | 0.00 | 0.00 | 0.00 | 0.00 |
| **7a** | 0.00 | 0.00 | 0.00 | 0.00 |
| **7b** | 0.00 | 0.00 | 0.00 | 0.00 |
| **7c** | -0.06 | -0.07 | 0.00 | 0.00 |
| **8a** | 0.00 | 0.00 | 0.00 | 0.00 |
| **8b** | 0.00 | 0.00 | 0.00 | 0.00 |
| **8c** | 0.00 | 0.00 | 0.00 | 0.00 |
| **9a** | 0.00 | 0.00 | 0.00 | 0.00 |
| **9b** | 0.00 | 0.00 | 0.00 | 0.00 |
| **9c** | 0.00 | 0.00 | 0.00 | 0.00 |
| **9d** | 0.00 | 0.00 | 0.00 | 0.00 |
| **10a** | 0.00 | 0.00 | 0.00 | 0.00 |
| **10b** | 0.00 | 0.00 | 0.00 | 0.00 |
| **10c** | 0.00 | 0.00 | 0.00 | 0.00 |
| **10d** | 0.00 | 0.00 | 0.00 | 0.00 |
| **Avg** | **-0.01** | **-0.04** | **0.00** | **0.00** |
| Max | 0.00 | 0.00 | 0.00 | 0.00 |
| Min | -0.15 | -1.08 | 0.00 | 0.00 |

| | Leftmost solution | | Rightmost solution | |
|---|---|---|---|---|
| **Val** | $E^1_{\sigma_H(x)}$ | $E^1_{\sigma_M(x)}$ | $E^2_{\sigma_H(x)}$ | $E^2_{\sigma_M(x)}$ |
| **1a** | 0.00 | 0.00 | 0.00 | 0.00 |
| **1b** | 0.00 | 0.00 | 0.00 | 0.00 |
| **1c** | 0.13 | 0.13 | 0.00 | 0.00 |
| **2a** | 0.00 | 0.00 | 0.00 | 0.00 |
| **2b** | 0.00 | 0.00 | 0.00 | 0.00 |
| **2c** | -0.02 | 0.00 | 0.59 | 0.00 |
| **3a** | 0.00 | 0.00 | 0.00 | 0.00 |
| **3b** | 0.01 | 0.00 | 0.03 | 0.00 |
| **3c** | 0.01 | 0.01 | 0.00 | 0.00 |
| **4a** | 0.05 | 0.00 | 0.07 | 0.00 |
| **4b** | 0.00 | 0.00 | 0.00 | 0.00 |
| **4c** | 0.00 | 0.00 | 0.00 | 0.00 |
| **4d** | -0.82 | 0.00 | -6.11 | 0.00 |
| **5a** | 0.00 | 0.00 | 0.00 | 0.00 |
| **5b** | 0.00 | 0.00 | 0.00 | 0.00 |
| **5c** | 0.00 | 0.00 | 0.00 | 0.00 |
| **5d** | 0.00 | 0.00 | 0.00 | 0.00 |
| **6a** | 0.00 | 0.00 | 0.00 | 0.00 |
| **6b** | 0.00 | 0.00 | 0.00 | 0.00 |
| **6c** | 0.14 | 0.13 | 0.00 | 0.00 |
| **7a** | 0.00 | 0.00 | 0.00 | 0.00 |
| **7b** | 0.00 | 0.00 | 0.00 | 0.00 |
| **7c** | -0.11 | 0.00 | 0.07 | 0.00 |
| **8a** | 0.00 | 0.00 | 0.00 | 0.00 |
| **8b** | 0.00 | 0.00 | 0.00 | 0.00 |
| **8c** | 0.05 | 0.00 | 0.07 | 0.00 |
| **9a** | 0.00 | 0.00 | 0.00 | 0.00 |
| **9b** | 0.00 | 0.00 | 0.00 | 0.00 |
| **9c** | 0.00 | 0.00 | 0.00 | 0.00 |
| **9d** | 0.03 | 0.00 | 0.26 | 0.00 |
| **10a** | 0.00 | 0.00 | 0.00 | 0.00 |
| **10b** | 0.00 | 0.00 | 0.00 | 0.00 |
| **10c** | 0.00 | 0.00 | 0.00 | 0.00 |
| **10d** | 0.00 | 0.00 | 0.00 | 0.00 |
| **Avg** | **-0.02** | **0.01** | **-0.15** | **0.00** |
| Max | 0.14 | 0.13 | 0.59 | 0.00 |
| Min | -0.11 | 0.00 | -6.11 | 0.00 |

## 3.3. CONCLUDING REMARKS ON QUALITY OF THE MATHEMATICAL FORMULAE

Table 5 sums up all the error in criteria evaluation using mathematical formulae. One can note the high quality of evaluation: the error is never greater than 0.3%. On average, considering both leftmost and rightmost solution, over the 4 criteria, the average error is about 0.06%. These results push into considering that mathematical formulae provide us high quality evaluation of criteria.

*Table 5. Quality of the mathematical formulas*

| Means of the instances | $E^1_{\overline{H(x)}}$ | $E^1_{\overline{M(x)}}$ | $E^2_{\overline{H(x)}}$ | $E^2_{\overline{M(x)}}$ | $E^1_{\sigma_H(x)}$ | $E^1_{\sigma_M(x)}$ | $E^2_{\sigma_H(x)}$ | $E^2_{\sigma_M(x)}$ |
|---|---|---|---|---|---|---|---|---|
| **Gdb** | 0.00 | 0.00 | 0.00 | 0.00 | 0.01 | 0.01 | 0.00 | 0.00 |
| **Val** | -0.01 | -0.04 | 0.00 | 0.00 | -0.02 | 0.01 | -0.15 | 0.00 |
| **Eglese** | 0.02 | 0.03 | -0.13 | -0.45 | 0.02 | -0.03 | -0.01 | -0.68 |

3.4. QUALITY OF THE SOLUTIONS OF THE GDB INSTANCES

Table 6 provides solutions for the 23 Gdb instances. When the obtained values are better or equal to solutions found solving the CARP by a bi-objective scheme, they are underlined. The rightmost solution solving the SCARP by a bi-objective approach is compared with the rightmost solution solving the CARP by a bi-objective approach. Similar comments remain true for the leftmost solutions. But, in the bi-objective resolution of the CARP, the standard deviations were not relevant. the comparison is achieved only on the expected cost and on the expected duration. The underline values report values better than (or equal to) the values obtained solving the SCARP by a mono-objective approach. The last column is the computational time in seconds. Note that the leftmost solutions found is 20 times better or equal using the SCARP bi-objective resolution than the CARP bi-objective resolution for the expected duration. The rightmost solutions are 15 times better than the best solutions found for the CARP bi-objective resolution in the expected cost. The SCARP resolution with a bi-objective resolution 9 times a best solution than the CARP resolution for the rightmost solution.

*Table 6. Solutions for the Gdb instances*

| | CARP resolution | | | | SCARP Resolution | | | | | |
|---|---|---|---|---|---|---|---|---|---|---|
| | Bi-objective resolution | | | | Mono-objective | Bi-objective resolution | | | | |
| **Gdb** | $h(x)^1$ | $m(x)^1$ | $h(x)^2$ | $m(x)^2$ | $\overline{H(x)}^{Mono}$ | $\overline{H(x)}^1$ | $\overline{M(x)}^1$ | $\overline{H(x)}^2$ | $\overline{M(x)}^2$ | $t$ |
| **1** | 316 | 74 | 337 | 63 | 337.0 | 337.0 | **63.0** | **337.0** | **63.0** | 52s |
| **2** | 339 | 69 | 395 | 59 | 388.0 | 367.0 | **62.0** | 417.0 | **59.0** | 58s |
| **3** | 275 | 65 | 339 | 59 | 296.0 | 296.0 | **60.0** | 347.0 | **59.0** | 52s |
| **4** | 287 | 74 | 350 | 64 | 313.0 | 317.0 | **72.0** | **350.0** | **64.0** | 46s |
| **5** | 377 | 78 | 447 | 64 | 409.0 | 433.0 | **68.0** | 479.0 | **64.0** | 58s |
| **6** | 298 | 75 | 351 | 64 | 324.0 | 324.0 | **68.0** | **351.0** | **64.0** | 52s |
| **7** | 325 | 68 | 381 | 61 | 351.0 | 359.0 | **68.0** | **368.0** | 66.0 | 58s |
| **8** | 350 | 44 | 390 | 38 | 372.1 | 399.0 | **44.0** | 493.0 | 39.0 | 92s |
| **9** | 309 | 43 | 333 | 37 | 326.3 | 361.2 | **38.0** | 379.0 | **37.0** | 100s |
| **10** | 275 | 71 | 297 | 54 | 283.0 | 283.0 | **66.0** | 329.0 | 57.0 | 58s |
| **11** | 395 | 81 | 421 | 64 | 396.0 | 409.0 | **81.0** | 439.0 | **61.0** | 92s |
| **12** | 458 | 97 | 547 | 91 | 534.0 | 523.7 | 99.1 | 613.0 | 93.0 | 56s |
| **13** | 544 | 128 | 544 | 128 | 552.0 | 556.0 | **128.0** | 556.0 | **128.0** | 70s |
| **14** | 100 | 21 | 112 | 17 | 96.6 | 104.0 | **20.0** | **108.0** | **17.0** | 52s |
| **15** | 58 | 15 | 60 | 13 | 58.0 | **58.0** | **15.0** | **60.0** | **13.0** | 52s |
| **16** | 127 | 27 | 135 | 19 | 129.0 | 131.0 | **20.0** | **135.0** | **18.0** | 62s |
| **17** | 91 | 15 | 91 | 15 | 91.0 | **91.0** | 17.0 | **91.0** | 17.0 | 66s |
| **18** | 164 | 33 | 178 | 27 | 161.5 | 172.0 | **33.0** | 182.0 | **27.0** | 84s |
| **19** | 55 | 21 | 63 | 17 | 63.0 | 63.0 | **17.0** | **63.0** | 18.0 | 36s |
| **20** | 121 | 36 | 131 | 20 | 123.0 | 123.0 | **29.0** | **127.0** | 21.0 | 52s |
| **21** | 156 | 30 | 160 | 22 | 154.6 | 160.0 | **28.0** | 166.0 | **22.0** | 76s |
| **22** | 200 | 26 | 207 | 20 | 201.1 | 207.0 | **26.0** | 208.0 | 20.0 | 94s |
| **23** | 235 | 23 | 241 | 20 | 237.1 | 239.0 | 28.0 | 245.0 | **19.0** | 120s |

Table 7 gives the 5 quality criteria for the solutions of the Gdb instances. For the leftmost solutions there is a gap of:

6.82% between the expected cost $\overline{H(x)}$ and the cost $h(x)$;

-4.99% between the expected duration $\overline{M(x)}$ and the duration $m(x)$;

2.04% between the expected cost $\overline{H(x)}$ and $\overline{H(x)}^{MONO}$.

Remember that $h(x)$ and $m(x)$ have been obtained by a bi-criteria scheme applied to the CARP. For the rightmost solution there is a gap of:

7.74% between the expected cost $\overline{H(x)}$ and the cost $h(x)$;

5.05 % between the expected duration $\overline{M(x)}$ and the duration $m(x)$.

*Table 7. Quality criteria of the solutions of Gdb instances (SCARP by a bi-objective scheme).*

| **Gdb** | $E^1_{\overline{H(x)}^1 - h(x)^1}$ | $E^1_{\overline{M(x)}^1 - m(x)^1}$ | $E^2_{\overline{H(x)}^2 - h(x)^2}$ | $E^2_{\overline{M(x)}^2 - m(x)^2}$ | $E^1_{\overline{H(x)}^1 - \overline{H(x)}^{MONO}}$ |
|---|---|---|---|---|---|
| **1** | 6.65 | -14.86 | 0.00 | 0.00 | 0.00 |
| **2** | 8.26 | -10.14 | 5.57 | 0.00 | -5.41 |
| **3** | 7.64 | -7.69 | 2.36 | 0.00 | 0.00 |
| **4** | 10.45 | -2.70 | 0.00 | 0.00 | 1.28 |
| **5** | 14.85 | -12.82 | 7.16 | 0.00 | 5.87 |
| **6** | 8.72 | -9.33 | 0.00 | 0.00 | 0.00 |
| **7** | 10.46 | 0.00 | -3.41 | 8.20 | 2.28 |
| **8** | 14.00 | 0.00 | 26.41 | 2.63 | 7.23 |
| **9** | 16.89 | -11.63 | 13.81 | 0.00 | 10.70 |
| **10** | 2.91 | -7.04 | 10.77 | 5.56 | 0.00 |
| **11** | 3.54 | 0.00 | 4.28 | -4.69 | 3.28 |
| **12** | 14.34 | 2.16 | 12.07 | 2.20 | -1.93 |
| **13** | 2.21 | 0.00 | 2.21 | 0.00 | 0.72 |
| **14** | 4.00 | -4.76 | -3.57 | 0.00 | 7.66 |
| **15** | 0.00 | 0.00 | 0.00 | 0.00 | 0.00 |
| **16** | 3.15 | -25.93 | 0.00 | -5.26 | 1.55 |
| **17** | 0.00 | 13.33 | 0.00 | 13.33 | 0.00 |
| **18** | 4.88 | 0.00 | 2.25 | 0.00 | 6.50 |
| **19** | 14.55 | -19.05 | 0.00 | 5.80 | 0.00 |
| **20** | 1.65 | -19.44 | -3.05 | 5.00 | 0.00 |
| **21** | 2.56 | -6.67 | 3.75 | 0.00 | 3.49 |
| **22** | 3.50 | 0.00 | 0.48 | 0.00 | 2.93 |
| **23** | 1.70 | 21.74 | 1.66 | -5.00 | 0.80 |
| **Avg** | **6.82** | **-4.99** | **3.60** | **1.21** | **2.04** |
| Max | 16.89 | 21.74 | 26.41 | 13.33 | 10.70 |
| Min | 0.00 | -19.44 | -3.57 | -5.00 | -5.41 |

3.5. QUALITY OF THE SOLUTIONS OF THE VAL INSTANCES

Table A2 (see Appendix) provides the results for the Val instances. One can note the high quality of the results obtained in front of those obtained solving the CARP with a bi-objective approach. The leftmost solutions are 4 times equal or better than solutions of the CARP using a bi-objective approach for the solution cost and 17 times better or equal for the solution duration of the longest trip. For the rightmost solutions. solution cost is 12 times better (or equal) solving the SCARP using a bi-objective approach than solving the CARP using a bi-objective one.

*Table 8. Solutions quality criteria of Val instances solving the SCARP by a bi-objective scheme*

| Val | $E^1_{\overline{H(x)}^1 - h(x)^1}$ | $E^1_{\overline{M(x)}^1 - m(x)^1}$ | $E^2_{\overline{H(x)}^2 - h(x)^2}$ | $E^2_{\overline{M(x)}^2 - m(x)^2}$ | $E^1_{\overline{H(x)}^1 - \overline{H(x)}^{MONO}}$ |
|---|---|---|---|---|---|
| **1a** | 0.00 | 1.72 | 0.00 | 1.72 | 0.00 |
| **1b** | 3.47 | -11.48 | -5.39 | 2.38 | 0.00 |
| **1c** | 15.92 | 2.44 | 14.52 | 5.00 | 5.97 |
| **2a** | 0.00 | 0.00 | -0.74 | 1.11 | 0.00 |
| **2b** | 3.08 | -9.90 | -0.33 | 2.56 | 3.08 |
| **2c** | 21.60 | 0.00 | 31.10 | 0.00 | -0.88 |
| **3a** | 0.00 | 0.00 | 0.00 | 0.00 | 0.00 |
| **3b** | 6.90 | 0.00 | 17.14 | 0.00 | 6.90 |
| **3c** | 27.54 | 0.00 | 27.54 | 0.00 | 8.64 |
| **4a** | 0.50 | 2.99 | 3.59 | 4.35 | 0.50 |
| **4b** | 7.77 | -12.38 | -0.43 | -1.20 | 5.21 |
| **4c** | 11.16 | 1.01 | 6.85 | 6.25 | 5.29 |
| **4d** | 12.10 | 2.62 | 18.74 | 0.00 | 4.82 |
| **5a** | 5.44 | -12.77 | 0.00 | 2.08 | 4.21 |
| **5b** | 6.28 | -13.39 | 0.00 | 3.49 | 5.29 |
| **5c** | 5.27 | 7.29 | -7.76 | 28.75 | 0.81 |
| **5d** | 2.35 | -2.47 | -11.22 | 9.72 | -3.50 |
| **6a** | 1.79 | 2.67 | -3.86 | 0.00 | 1.79 |
| **6b** | 5.15 | -20.59 | 2.66 | 2.00 | 0.00 |
| **6c** | 13.25 | -18.18 | 9.73 | 0.00 | 3.43 |
| **7a** | 1.43 | -30.59 | -2.08 | 0.00 | 1.43 |
| **7b** | 0.00 | 6.90 | 4.01 | 7.84 | 0.00 |
| **7c** | 7.16 | -10.00 | 22.73 | -2.50 | 1.99 |
| **8a** | 2.33 | -21.71 | 0.47 | 1.15 | 0.74 |
| **8b** | 6.08 | -3.00 | 2.42 | 2.53 | 2.44 |
| **8c** | 12.11 | -4.05 | 8.69 | 0.00 | 4.80 |
| **9a** | 1.53 | 8.54 | 4.20 | 2.94 | 2.48 |
| **9b** | 3.37 | -14.63 | 8.24 | 6.90 | 1.51 |
| **9c** | 5.12 | -5.80 | 0.00 | 1.96 | 2.65 |
| **9d** | 10.53 | -6.00 | 8.29 | 0.00 | 4.75 |
| **10a** | 2.80 | -7.69 | 3.79 | 6.59 | 2.80 |
| **10b** | 4.36 | -13.51 | 5.01 | 9.09 | 4.36 |
| **10c** | 4.46 | -3.23 | 5.62 | 9.09 | 1.96 |
| **10d** | 12.29 | 6.56 | 6.22 | 5.56 | 7.30 |
| **Avg** | **6.56** | **-5.25** | **5.29** | **3.51** | **2.67** |
| Max | 27.54 | 7.29 | 31.10 | 28.75 | 8.64 |
| Min | 0.00 | -30.59 | -11.22 | -2.50 | -3.50 |

Table 8 gives the solutions quality criteria for the Gdb instances. For the leftmost solution. there is a gap of: 6.56% between the expected cost $\overline{H(x)}$ and the cost $h(x)$; -5.25% between the expected duration $\overline{M(x)}$ and the duration $m(x)$; 2.67% between the expected cost $\overline{H(x)}$ and $\overline{H(x)}^{MONO}$. For the rightmost solution. there is a gap of: 5.29% between the expected cost $\overline{H(x)}$ and the cost $h(x)$; 3.51 % between the expected duration $\overline{M(x)}$ and the duration $m(x)$.

### 3.6. QUALITY OF SOLUTIONS FOR THE EGLESE INSTANCES

Results for the Eglese instances are provided in table A3 and table 9 provides results on the gap between the two criteria for both the leftmost and the rightmost solution between the bi-objective resolution of the CARP and of the SCARP. The gap from solution duration is not relevant: around 3% for both leftmost and rightmost solutions between the bi-objective resolution of the CARP and the bi-objective resolution of the SCARP. For the solution cost one can note a gap around 30% for both leftmost and rightmost solutions but gap of only 6.44% with the best solution found solving the SCARP with a mono-objective approach.

*Table 9. Solutions quality criteria of Val instances solving the SCARP by a bi-objective scheme*

| egl- | $E^1_{\overline{H(x)}^1 - h(x)^1}$ | $E^1_{\overline{M(x)}^1 - m(x)^1}$ | $E^2_{\overline{H(x)}^2 - h(x)^2}$ | $E^2_{\overline{M(x)}^2 - m(x)^2}$ | $E^1_{\overline{H(x)}^1 - \overline{H(x)}^{MONO}}$ |
|---|---|---|---|---|---|
| e1-A | 17.90 | -13.04 | 9.39 | 0.00 | 3.64 |
| e1-B | 25.50 | -2.26 | 24.19 | 0.00 | 15.13 |
| e1-C | 28.67 | 0.00 | 32.23 | 0.00 | 4.08 |
| e2-A | 23.22 | -8.71 | 27.09 | 3.66 | 9.84 |
| e2-B | 24.30 | 47.34 | 26.85 | 0.00 | 3.23 |
| e2-C | 32.32 | 2.81 | 33.11 | 4.16 | 6.01 |
| e3-A | 23.37 | -4.47 | 5.12 | 0.85 | 10.88 |
| e3-B | 27.84 | 0.24 | 28.86 | 3.41 | -1.47 |
| e3-C | 28.83 | -2.87 | 39.79 | 2.32 | 10.49 |
| e4-A | 27.76 | 4.04 | 21.80 | 2.32 | 7.86 |
| e4-B | 29.32 | 0.46 | 30.92 | 6.34 | 5.56 |
| e4-C | 31.77 | 7.16 | 42.81 | 5.73 | 0.80 |
| s1-A | 25.66 | 3.91 | 14.18 | 0.00 | 1.44 |
| s1-B | 27.89 | -3.86 | 21.84 | 0.00 | 2.46 |
| s1-C | 23.42 | 4.65 | 42.61 | 2.85 | -2.07 |
| s2-A | 32.10 | -1.70 | 12.59 | 4.90 | 18.37 |
| s2-B | 34.57 | -4.73 | 48.49 | 0.00 | 4.76 |
| s2-C | 34.07 | 5.36 | 45.42 | 0.92 | 5.55 |
| s3-A | 32.35 | -1.46 | 19.28 | 5.52 | 14.20 |
| s3-B | 33.86 | 0.10 | 36.84 | 0.00 | 0.15 |
| s3-C | 34.29 | 6.23 | 40.38 | 0.00 | 4.77 |
| s4-A | 29.67 | 21.16 | 32.34 | 1.27 | 14.37 |
| s4-B | 34.90 | 4.51 | 48.13 | 4.67 | 5.34 |
| s4-C | 39.46 | 13.51 | 52.56 | 24.59 | 9.14 |
| **Avg** | **29.29** | **3.27** | **30.70** | **3.06** | **6.44** |
| **Max** | 39.46 | 47.34 | 52.56 | 24.59 | 15.13 |
| **Min** | 17.90 | -13.04 | 5.12 | 0.00 | 0.15 |

### 3.7. CONCLUDING REMARKS ON SOLUTIONS QUALITY CRITERIA

Table 10 provides the average gap over the 5 criteria for all the instances including Gdb. Val and Eglese instances. On average. one can note that for the leftmost solution. there is a gap of:

14.22% between the expected cost $\overline{H(x)}$ and the cost $h(x)$;

-2.32% between the expected duration $\overline{M(x)}$ and the duration $m(x)$;

5.31% between the expected cost $\overline{H(x)}$ and $\overline{H(x)}^{MONO}$;

Remember that $h(x)$ and $m(x)$ have been obtained by a bi-criteria scheme applied to the CARP. For the rightmost solution. there is a gap of:

14.58% between the expected cost $\overline{H(x)}$ and the cost $h(x)$;

3.87 % between the expected duration $\overline{M(x)}$ and the duration $m(x)$

*Table 10. Solutions quality criteria of Val instances solving the SCARP by a bi-objective scheme*

| Instances | $E^1_{\overline{H(x)}^1 - h(x)^1}$ | $E^1_{\overline{M(x)}^1 - m(x)^1}$ | $E^2_{\overline{H(x)}^2 - h(x)^2}$ | $E^2_{\overline{M(x)}^2 - m(x)^2}$ | $E^1_{\overline{H(x)}^1 - \overline{H(x)}^{MONO}}$ |
|---|---|---|---|---|---|
| **Gdb** | 6.82 | -4.99 | 3.60 | 1.21 | 2.04 |
| **Val** | 6.56 | -5.25 | 5.29 | 3.51 | 2.67 |
| **Eglese** | 29.29 | 3.27 | 30.70 | 3.06 | 6.44 |

Let us note that Eglese instances seem to be more difficult to solve since a gap of around 30% is expected for both $E^1_{\overline{H(x)}^1 - h(x)^1}$ and $E^2_{\overline{H(x)}^2 - h(x)^2}$. These value must be analyzed taking into account the short gap of 6% between the mono-objective resolution of the SCARP and the bi-objective resolution one. These results push into accepting that random events consequences are more dramatic for Eglese's instances which are true CARP with non required arcs.

## 4. CONCLUDING REMARKS

It is possible to solve a multi-objective and stochastic CARP problem as soon as the laws of the random variables modeling the demand are known. The required hypotheses (a trip requires at most a move to the depot, and this move is just before the last arc to collect) permit to introduce mathematical analysis to obtain satisfactory estimation of criteria to minimize. The formalization concerns a Gaussian law since the demand on each arc can be considered as the sum of independent random variables (having a mean and a standard deviation).

The NSGA-II template is used to simultaneously optimize both solution cost and the length of the longest trip. The experiments were carried out using the well known standard benchmarks of the state of the art literature on CARP.

The results prove that mathematical formulae are high quality ones and that the NSGA II template is able to optimize the SCARP : the stochastic solution of the first front are very close to the best one found solving the CARP by the NSGA II template and very close to the best solution found solving the SCARP by a mono-objective approach.

The work presented is a step forward stochastic resolution of routing problem with the aim to obtain in rather short computation time, robust solution on several criteria. Our work is now directed to:

- extension of the NSGA-II template for optimization of more than two criteria simultaneously;
- extension of the mathematical formulae to address the problem of a heterogonous fleet of vehicles.

## 6. APPENDIX

*Table A1. Standard deviation of cost and duration for the leftmost and rightmost solutions according to mathematical expressions and replications values for the Val files.*

| | Values obtained by mathematical expression | | | | Values obtained by replications | | | |
|---|---|---|---|---|---|---|---|---|
| | Standard deviation cost | | Standard deviation duration | | Standard deviation cost | | Standard deviation duration | |
| Val | $\sigma_H^1(x)$ | $\sigma_H^2(x)$ | $\sigma_M^1(x)$ | $\sigma_M^2(x)$ | $\sigma_H^1(x,n)$ | $\sigma_H^2(x,n)$ | $\sigma_M^1(x,n)$ | $\sigma_M^2(x,n)$ |
| 1a | 0.00 | 0.00 | 0.00 | 0.00 | 0.00 | 0.00 | 0.00 | 0.00 |
| 1b | 0.00 | 0.00 | 0.00 | 0.00 | 0.00 | 0.00 | 0.00 | 0.00 |
| 1c | 0.38 | 0.38 | 0.00 | 0.00 | 0.00 | 0.00 | 0.00 | 0.00 |
| 2a | 0.00 | 0.00 | 0.00 | 0.00 | 0.00 | 0.00 | 0.00 | 0.00 |
| 2b | 0.00 | 0.00 | 0.00 | 0.00 | 0.00 | 0.00 | 0.00 | 0.00 |
| 2c | 0.59 | 0.00 | 0.42 | 0.00 | 0.72 | 0.00 | 0.00 | 0.00 |
| 3a | 0.00 | 0.00 | 0.00 | 0.00 | 0.00 | 0.00 | 0.00 | 0.00 |
| 3b | 0.30 | 0.00 | 0.16 | 0.00 | 0.29 | 0.00 | 0.15 | 0.00 |
| 3c | 0.02 | 0.02 | 0.00 | 0.00 | 0.00 | 0.00 | 0.00 | 0.00 |
| 4a | 0.22 | 0.00 | 0.10 | 0.00 | 0.00 | 0.00 | 0.00 | 0.00 |
| 4b | 0.00 | 0.00 | 0.00 | 0.00 | 0.00 | 0.00 | 0.00 | 0.00 |
| 4c | 0.00 | 0.00 | 0.00 | 0.00 | 0.00 | 0.00 | 0.00 | 0.00 |
| 4d | 1.30 | 0.10 | 1.10 | 0.00 | 6.28 | 0.09 | 6.12 | 0.00 |
| 5a | 0.00 | 0.00 | 0.00 | 0.00 | 0.00 | 0.00 | 0.00 | 0.00 |
| 5b | 0.00 | 0.00 | 0.00 | 0.00 | 0.00 | 0.00 | 0.00 | 0.00 |
| 5c | 0.00 | 0.00 | 0.00 | 0.00 | 0.00 | 0.00 | 0.00 | 0.00 |
| 5d | 0.00 | 0.00 | 0.00 | 0.00 | 0.00 | 0.00 | 0.00 | 0.00 |
| 6a | 0.00 | 0.00 | 0.00 | 0.00 | 0.00 | 0.00 | 0.00 | 0.00 |
| 6b | 0.00 | 0.00 | 0.00 | 0.00 | 0.00 | 0.00 | 0.00 | 0.00 |
| 6c | 0.49 | 0.47 | 0.00 | 0.00 | 0.00 | 0.00 | 0.00 | 0.00 |
| 7a | 0.00 | 0.00 | 0.00 | 0.00 | 0.00 | 0.00 | 0.00 | 0.00 |
| 7b | 0.00 | 0.00 | 0.00 | 0.00 | 0.00 | 0.00 | 0.00 | 0.00 |
| 7c | 0.76 | 0.00 | 0.26 | 0.00 | 1.17 | 0.00 | 0.23 | 0.00 |
| 8a | 0.00 | 0.00 | 0.00 | 0.00 | 0.00 | 0.00 | 0.00 | 0.00 |
| 8b | 0.00 | 0.00 | 0.00 | 0.00 | 0.00 | 0.00 | 0.00 | 0.00 |
| 8c | 0.29 | 0.03 | 0.05 | 0.00 | 0.00 | 0.00 | 0.00 | 0.00 |
| 9a | 0.00 | 0.00 | 0.00 | 0.00 | 0.00 | 0.00 | 0.00 | 0.00 |
| 9b | 0.00 | 0.00 | 0.00 | 0.00 | 0.00 | 0.00 | 0.00 | 0.00 |
| 9c | 0.00 | 0.00 | 0.00 | 0.00 | 0.00 | 0.00 | 0.00 | 0.00 |
| 9d | 0.12 | 0.00 | 0.12 | 0.00 | 0.00 | 0.00 | 0.00 | 0.00 |
| 10a | 0.00 | 0.00 | 0.00 | 0.00 | 0.00 | 0.00 | 0.00 | 0.00 |
| 10b | 0.00 | 0.00 | 0.00 | 0.00 | 0.00 | 0.00 | 0.00 | 0.00 |
| 10c | 0.00 | 0.00 | 0.00 | 0.00 | 0.00 | 0.00 | 0.00 | 0.00 |
| 10d | 0.00 | 0.00 | 0.00 | 0.00 | 0.00 | 0.00 | 0.00 | 0.00 |

*Table A2. Solutions for the Val instances compared to previous published results*

| | CARP resolution | | | | SCARP Resolution | | | | | |
|---|---|---|---|---|---|---|---|---|---|---|
| | Bi-objective resolution | | | | Mono-objective | Bi-objective resolution | | | | |
| **Val** | $h(x)^1$ | $m(x)^1$ | $h(x)^2$ | $m(x)^2$ | $\overline{H(x,n)}^{Mono}$ | $\overline{H(x)}^1$ | $\overline{M(x)}^1$ | $\overline{H(x)}^2$ | $\overline{M(x)}^2$ | $t$ |
| **1a** | 173 | 58 | 173 | 58 | 173.0 | **173.0** | 59.0 | **173.0** | 59.0 | 97s |
| **1b** | 173 | 61 | 204 | 42 | 179.0 | 179.0 | **54.0** | 193.0 | 43.0 | 111s |
| **1c** | 245 | 41 | 248 | 40 | 268.0 | 284.0 | 42.0 | 284.0 | 42.0 | 91s |
| **2a** | 227 | 114 | 270 | 90 | 227.0 | **227.0** | **114.0** | **268.0** | 91.0 | 78s |
| **2b** | 260 | 101 | 307 | 78 | 260.0 | 268.0 | **91.0** | **306.0** | 80.0 | 75s |
| **2c** | 463 | 71 | 463 | 71 | 568.0 | 563.0 | **71.0** | 607.0 | **71.0** | 74s |
| **3a** | 81 | 41 | 88 | 31 | 81.0 | **81.0** | **41.0** | **88.0** | **31.0** | 80s |
| **3b** | 87 | 32 | 105 | 27 | 87.0 | 93.0 | 32.0 | 123.0 | **27.0** | 79s |
| **3c** | 138 | 27 | 138 | 27 | 162.0 | 176.0 | **27.0** | 176.0 | **27.0** | 92s |
| **4a** | 400 | 134 | 446 | 92 | 400.0 | 402.0 | 138.0 | 462.0 | 96.0 | 155s |
| **4b** | 412 | 105 | 468 | 83 | 422.0 | 444.0 | **92.0** | **466.0** | **82.0** | 152s |
| **4c** | 430 | 99 | 482 | 80 | 454.0 | 478.0 | 100.0 | 515.0 | 85.0 | 160s |
| **4d** | 539 | 80 | 539 | 80 | 576.4 | 604.2 | 82.1 | 640.0 | 80.0 | 158s |
| **5a** | 423 | 141 | 474 | 96 | 428.0 | 446.0 | **123.0** | **474.0** | 98.0 | 147s |
| **5b** | 446 | 112 | 506 | 86 | 450.2 | 474.0 | **97.0** | **506.0** | 89.0 | 155s |
| **5c** | 474 | 96 | 541 | 80 | 495.0 | 499.0 | 103.0 | **499.0** | 103.0 | 146s |
| **5d** | 595 | 81 | 686 | 72 | 631.1 | 609.0 | 79.0 | **609.0** | 79.0 | 131s |
| **6a** | 223 | 75 | 259 | 56 | 223.0 | 227.0 | 77.0 | **249.0** | **56.0** | 110s |
| **6b** | 233 | 68 | 263 | 50 | 245.0 | 245.0 | **54.0** | 270.0 | 51.0 | 110s |
| **6c** | 317 | 55 | 329 | 45 | 347.1 | 359.0 | **45.0** | 361.0 | **45.0** | 163s |
| **7a** | 279 | 85 | 289 | 59 | 279.0 | 283.0 | **59.0** | **283.0** | **59.0** | 167s |
| **7b** | 283 | 58 | 299 | 51 | 283.0 | **283.0** | 62.0 | 311.0 | 55.0 | 171s |
| **7c** | 335 | 50 | 352 | 40 | 352.0 | 359.0 | **45.0** | 432.0 | **39.0** | 150s |
| **8a** | 386 | 129 | 429 | 87 | 392.1 | 395.0 | **101.0** | 431.0 | 88.0 | 142s |
| **8b** | 395 | 100 | 455 | 79 | 409.0 | 419.0 | 97.0 | 466.0 | 81.0 | 135s |
| **8c** | 545 | 74 | 610 | 67 | 583.0 | 611.0 | **71.0** | 663.0 | **67.0** | 134s |
| **9a** | 326 | 82 | 333 | 68 | 323.0 | 331.0 | 89.0 | 347.0 | 70.0 | 230s |
| **9b** | 326 | 82 | 340 | 58 | 332.0 | 337.0 | **70.0** | 368.0 | 62.0 | 214s |
| **9c** | 332 | 69 | 389 | 51 | 340.0 | 349.0 | **65.0** | **389.0** | 52.0 | 221s |
| **9d** | 399 | 50 | 434 | 44 | 421.0 | 441.0 | **47.0** | 470.0 | **44.0** | 210s |
| **10a** | 428 | 143 | 449 | 91 | 428.0 | 440.0 | **132.0** | 466.0 | 97.0 | 236s |
| **10b** | 436 | 111 | 459 | 77 | 436.0 | 455.0 | **96.0** | 482.0 | 84.0 | 230s |
| **10c** | 448 | 93 | 498 | 66 | 459.0 | 468.0 | 90.0 | 526.0 | 72.0 | 223s |
| **10d** | 537 | 61 | 595 | 54 | 562.0 | 603.0 | 65.0 | 632.0 | 57.0 | 235s |

*Table A3. Solutions for the Eglese's instances compared to previous published results*

| | **CARP resolution** | | | | **SCARP Resolution** | | | | | |
|---|---|---|---|---|---|---|---|---|---|---|
| | **Bi-objective resolution** | | | | **Mono-objective** | **Bi-objective resolution** | | | | |
| **egl-** | $h(x)^1$ | $m(x)^1$ | $h(x)^2$ | $m(x)^2$ | $\overline{H(x,n)}^{Mono}$ | $\overline{H(x)}^1$ | $\overline{M(x)}^1$ | $\overline{H(x)}^2$ | $\overline{M(x)}^2$ | $t$ |
| **e1-A** | 3548 | 943 | 3824 | 820 | 4036.0 | 4183.0 | 820.0 | 4183.0 | **820.0** | 195s |
| **e1-B** | 4525 | 839 | 4573 | 820 | 4932.7 | 5679.0 | **820.0** | 5679.0 | **820.0** | 164s |
| **e1-C** | 5687 | 836 | 5764 | 820 | 7030.8 | 7317.6 | **836.0** | 7621.6 | **820.0** | 165s |
| **e2-A** | 5018 | 953 | 6072 | 820 | 5629.0 | 6183.0 | **870.0** | 7717.0 | **850.0** | 148s |
| **e2-B** | 6411 | 564 | 6810 | 820 | 7720.1 | 7969.1 | 831.0 | 8638.2 | 820.0 | 170s |
| **e2-C** | 8440 | 854 | 8651 | 820 | 10534.5 | 11168.0 | 878.0 | 11515.4 | 854.1 | 154s |
| **e3-A** | 5956 | 917 | 7935 | 820 | 6627.1 | 7348.0 | **876.0** | 8341.0 | 827.0 | 180s |
| **e3-B** | 7911 | 872 | 8455 | 820 | 10264.3 | 10113.2 | 874.1 | 10895.0 | 848.0 | 186s |
| **e3-C** | 10349 | 864 | 10511 | 820 | 12066.3 | 13332.5 | **839.2** | 14693.2 | 839.0 | 225s |
| **e4-A** | 6548 | 890 | 7362 | 820 | 7756.5 | 8366.0 | 926.0 | 8967.0 | 839.0 | 260s |
| **e4-B** | 9116 | 874 | 9584 | 820 | 11168.4 | 11789.0 | 878.0 | 12547.8 | 872.0 | 237s |
| **e4-C** | 11802 | 820 | 11802 | 820 | 15427.9 | 15551.4 | 878.7 | 16854.4 | 867.0 | 218s |
| **s1-A** | 5102 | 1023 | 6582 | 924 | 6320.5 | 6411.3 | 1063.0 | 7515.0 | **924.0** | 151s |
| **s1-B** | 6500 | 984 | 8117 | 912 | 8113.5 | 8313.1 | **946.0** | 9890.0 | **912.0** | 150s |
| **s1-C** | 8694 | 946 | 9205 | 912 | 10957.6 | 10730.5 | 990.0 | 13127.0 | 938.0 | 160s |
| **s2-A** | 10207 | 1058 | 12222 | 979 | 11390.7 | 13483.0 | **1040.0** | 13761.0 | 1027.0 | 315s |
| **s2-B** | 13548 | 1058 | 14334 | 979 | 17403.1 | 18232.0 | **1008.0** | 21284.8 | **979.0** | 300s |
| **s2-C** | 16932 | 1040 | 16975 | 979 | 21507.3 | 22701.5 | 1095.7 | 24684.3 | 988.0 | 352s |
| **s3-A** | 10456 | 1099 | 12605 | 979 | 12118.6 | 13838.9 | **1083.0** | 15035.2 | 1033.0 | 404s |
| **s3-B** | 14004 | 1040 | 15103 | 979 | 18717.6 | 18745.4 | 1041.0 | 20666.7 | 979.0 | 336s |
| **s3-C** | 17825 | 998 | 18043 | 979 | 22847.5 | 23936.6 | 1060.2 | 25329.4 | 979.0 | 406s |
| **s4-A** | 12730 | 1040 | 12912 | 1027 | 14433.8 | 16507.3 | 1260.1 | 17088.1 | 1040.0 | 397s |
| **s4-B** | 16792 | 1027 | 16792 | 1027 | 21503.6 | 22651.6 | 1073.3 | 24873.4 | 1075.0 | 411s |
| **s4-C** | 21309 | 1027 | 21309 | 1027 | 27226.9 | 29716.5 | 1165.7 | 32509.7 | 1279.5 | 474s |